\documentclass[11pt, reqno]{amsart}
\usepackage{mathrsfs}
\usepackage[active]{srcltx}
\usepackage{mathrsfs}
\usepackage{longtable}

\usepackage{color}
\usepackage[unicode]{hyperref}
\hypersetup{
    colorlinks = true,%
    citecolor = [rgb]{0.0,0.0,0.9},
    filecolor=black,%
    linkcolor = [rgb]{0.65,0.0,0.0},%
    anchorcolor = red,
    pagecolor = red,
    urlcolor= [rgb]{0.65,0.0,0.0}
}

\newtheorem{theorem}{Theorem}[section]
\newtheorem{lemma}{Lemma}[section]
\newtheorem{corollary}{Corollary}[section]

\newtheorem{remark}{Remark}[section]
\numberwithin{equation}{section}

\begin{document}

\title[Gagliardo-Nirenberg inequalities on metric measure
spaces]{Metric measure spaces supporting
Gagliardo-Nirenberg inequalities: volume non-collapsing and rigidities 
}

\author{Alexandru Krist\'aly}
\address{Department of Economics, Babe\c s-Bolyai University, 400591 Cluj-Napoca,
	Romania \&  Institute of Applied Mathematics, \'Obuda
University, 1034 
Budapest, Hungary}

\email{alex.kristaly@econ.ubbcluj.ro; alexandrukristaly@yahoo.com}



\begin{abstract}
{\footnotesize \noindent Let $({M},\textsf{d},\textsf{m})$ be a
metric measure space which satisfies the Lott-Sturm-Villani
cur\-va\-ture-dimension condition $\textsf{CD}(K,n)$ for some $K\geq
0$ and $n\geq 2$, and a lower $n-$density assumption  at some point
of $M$. We prove that if $({M},\textsf{d},\textsf{m})$ supports the
Gagliardo-Nirenberg inequality or any of its limit cases
($L^p-$logarithmic Sobolev inequality or Faber-Krahn-type
inequality), then a  {\it global non-collapsing $n-$dimensional
volume growth} holds, i.e., there exists a universal constant
$C_0>0$ such that
 $\textsf{m}( B_x(\rho))\geq C_0 \rho^n$  for all $x\in {M}$ and $\rho\geq 0,$
 where $B_x(\rho)=\{y\in M:{\sf d}(x,y)<\rho\}$. Due to the quantitative character of the volume growth estimate,
 we establish several rigidity results on
Riemannian manifolds with non-negative Ricci curvature supporting
Gagliardo-Nirenberg
  inequalities by exploring a quantitative  Perelman-type homotopy construction developed by  Munn (J. Geom. Anal., 2010).
Further rigidity results are also presented on some reversible
Finsler manifolds.

 }
\end{abstract}
\vspace*{-2cm} \maketitle

\tableofcontents


\vspace*{-1.0cm}
 \section{\sc Introduction}

An important role in the theory of geometric functional inequalities
is played by the Gagliardo-Nirenberg interpolation inequality and
its limit cases. The present paper is devoted to the study of
Gagliardo-Nirenberg inequalities on metric measure spaces; to be
more precise, we shall
\begin{itemize}
  \item[(a)] establish {\it quantitative volume non-collapsing properties} of  metric measure spaces satisfying the Lott-Sturm-Villani curvature-dimension
  condition  $\textsf{CD}(K,n)$ for some $K\geq 0$ and $n\geq 2,$
   in the presence of a Gagliardo-Nirenberg inequality or one of its limit cases ($L^p-$logarithmic Sobolev inequality or Faber-Krahn-type inequality);
  \item[(b)] provide {\it rigidity} results in the framework of  Riemannian and Finsler  manifolds with non-negative Ricci curvature which support {\it $($almost$)$ optimal} Gagliardo-Nirenberg inequalities by using the volume non-collapsing property  from
  (a) and a quantitative homotopy  construction due to Munn \cite{Munn-JGA} and Perelman \cite{Perelman-JAMS}.
\end{itemize}
In \S \ref{sect1.1} we recall the optimal Gagliardo-Nirenberg
inequalities on normed spaces which play a comparison role in our
 investigations; in \S \ref{sect1.2}  we
 present the main results of the paper.

\subsection{Recalling optimal Gagliardo-Nirenberg inequalities on normed
spaces}\label{sect1.1} The optimal Gagli\-ardo-Nirenberg inequality
in the Euclidean case has been obtained by Del Pino and Dolbeault
\cite{DelPino-Dolb} for a certain range of parameters by using
symmetrization arguments. By using  mass transportation argument,
Cordero-Erausquin, Nazaret and Villani \cite{CE-N-Villani} extended
the results from \cite{DelPino-Dolb} to prove optimal
Gagliardo-Nirenberg inequalities on arbitrary normed spaces. In the
sequel, we recall the main theorems from \cite{CE-N-Villani} and
some related results.

 Let $\|\cdot\|$ be an arbitrary norm on $\mathbb R^n;$ without loss of generality, we may assume that  the
Lebesgue measure of the unit ball in $(\mathbb R^n,\|\cdot\|)$ is
the volume of the $n-$di\-men\-sional Euclidean unit ball
$\omega_n={\pi^\frac{n}{2}}{\Gamma(\frac{n}{2}+1)}^{-1}$. The  dual
norm $\|\cdot\|_*$ of $\|\cdot\|$ is given by
$\|x\|_*=\sup_{\|y\|\leq 1}x\cdot y$ where $'\cdot'$ is the
Euclidean inner product. Let $p\in [1,n)$ and $L^p(\mathbb R^n)$ be
the Lebesgue space of order $p$. As usual, we consider the Sobolev
spaces
$$\dot W^{1,p}(\mathbb R^n)=\{u\in L^{p^\star}(\mathbb R^n):\nabla u\in L^p(\mathbb
R^n)\}$$ and $$W^{1,p}(\mathbb R^n)=\{u\in L^{p}(\mathbb R^n):\nabla
u\in L^p(\mathbb R^n)\},$$  where $p^\star=\frac{pn}{n-p}$ and
$\nabla$ is the gradient operator. On account of the Finslerian
duality (see also \S \ref{sect-finsler}), if $u\in \dot
W^{1,p}(\mathbb R^n),$ the norm of $\nabla u$ is defined by
$$\|\nabla u\|_{L^p}=\left(\int_{\mathbb R^n}\|\nabla u(x)\|_*^pdx\right)^{1/p},$$
where $dx$ is the Lebesgue measure on $\mathbb R^n$.

Fix $n\geq 2,$ $p\in (1,n)$ and $\alpha\in
(0,\frac{n}{n-p}]\setminus \{1\}$; for every $\lambda>0,$ let
$$h_{\alpha,p}^\lambda(x)=(\lambda+(\alpha-1)\|x\|^{p'})_+^\frac{1}{1-\alpha},\ x\in \mathbb
R^n,\footnote{The  function $h_{\alpha,p}^\lambda$ is positive
everywhere for $\alpha>1$ while $h_{\alpha,p}^\lambda$ has always a
compact support for $\alpha<1$.}$$ where $p'=\frac{p}{p-1}$ is the
conjugate to $p$, and $r_+=\max\{0,r\}$ for $r\in \mathbb R.$
The following {\it optimal Gagliardo-Nirenberg inequalities} are
known on normed spaces:
\\

\noindent {\bf Theorem A.} {\rm [see  \cite[Theorem
4]{CE-N-Villani}]} {\it Let $n\geq 2$, $p\in (1,n)$ and $\|\cdot\|$
be an arbitrary norm on $\mathbb R^n$.}
\begin{itemize}
 {\it  \item If $1<\alpha\leq \frac{n}{n-p}$, then
  \begin{equation}\label{Villani-1}
    \|u\|_{L^{\alpha p}}\leq \mathcal G_{\alpha,p,n}
    \|\nabla u\|_{L^p}^{\theta}\|u\|_{L^{\alpha(p-1)+1}}^{1-\theta},\
    \forall u\in \dot W^{1,p}(\mathbb R^n),
  \end{equation}
 where
 \begin{equation}\label{theta-best}
    \theta=\frac{p^\star(\alpha-1)}{\alpha p(p^\star-\alpha
 p+\alpha-1)},
 \end{equation}
  and the best constant $$\mathcal
 G_{\alpha,p,n}=\left(\frac{\alpha-1}{p'}\right)^\theta
 \frac{\left(\frac{p'}{n}\right)^{\frac{\theta}{p}+\frac{\theta}{n}}\left(\frac{\alpha (p-1)+1}{\alpha
 -1}-\frac{n}{p'}\right)^\frac{1}{\alpha p}
 \left(\frac{\alpha (p-1)+1}{\alpha
 -1}\right)^{\frac{\theta}{p}-\frac{1}{\alpha p}}}{\left(\omega_n {\sf B}\left(\frac{\alpha (p-1)+1}{\alpha
 -1}-\frac{n}{p'},\frac{n}{p'}\right)\right)^{\frac{\theta}{n}}}$$
 is achieved by the family of functions
 $h_{\alpha,p}^\lambda$,
 $\lambda>0;$
  \item If $0<\alpha<1$, then
  \begin{equation}\label{Villani-2}
    \|u\|_{L^{\alpha(p-1)+1}}\leq \mathcal N_{\alpha,p,n}
    \|\nabla  u\|_{L^p}^{\gamma}\|u\|_{L^{\alpha p}}^{1-\gamma},\
    \forall u\in \dot W^{1,p}(\mathbb R^n),
  \end{equation}
 where
\begin{equation}\label{gamma-best}
     \gamma=\frac{p^\star(1-\alpha)}{(p^\star-\alpha
 p)(\alpha p+1-\alpha)},
\end{equation}
  and the best constant $$\mathcal
N_{\alpha,p,n}={\small \left(\frac{1-\alpha}{p'}\right)^\gamma
 \frac{\left(\frac{p'}{n}\right)^{\frac{\gamma}{p}+\frac{\gamma}{n}}\left(\frac{\alpha (p-1)+1}{1-\alpha
 }+\frac{n}{p'}\right)^{\frac{\gamma}{ p}-\frac{1}{\alpha(p-1)+1}}
 \left(\frac{\alpha (p-1)+1}{1-\alpha
 }\right)^{\frac{1}{\alpha(p-1)+1}}}{\left(\omega_n {\sf B}\left(\frac{\alpha (p-1)+1}{1-\alpha
 },\frac{n}{p'}\right)\right)^{\frac{\gamma}{n}}}}$$
 is achieved  by the family of functions $h_{\alpha,p}^\lambda$,
 $\lambda>0.$
}
\end{itemize}
Hereafter, ${\sf B}(\cdot,\cdot)$ is the Euler
beta-function.\newpage

The borderline case  $\alpha=\frac{n}{n-p}$ (thus $\theta=1$)
reduces to the {\it optimal Sobolev inequality}, see Aubin
\cite{Aubin} and Talenti \cite{Talenti} in the Euclidean case, and
Alvino, Ferone, Lions and  Trombetti \cite{AIHP-Lions} for normed
spaces. Furthermore, inequalities (\ref{Villani-1}) and
(\ref{Villani-2}) degenerate to the {\it optimal $L^p-$logarithmic
Sobolev inequality} whenever $\alpha\to 1$ (called also as the
entropy-energy inequality involving the Shannon entropy), while
(\ref{Villani-2}) reduces to a {\it Faber-Krahn-type inequality}
whenever $\alpha\to 0$,
respectively. More precisely, one has\\

\noindent {\bf Theorem B.} {\it Let $n\geq 2$, $p\in (1,n)$ and
$\|\cdot\|$ be an arbitrary  norm on $\mathbb R^n$.}
\begin{itemize}
  \item {\bf Limit case I} ($\alpha\to 1$) {\rm [see  \cite[Theorem 1.1]{Gentil}\footnote{Gentil \cite{Gentil} proved an optimal $L^p-$logarith\-mic So\-bo\-lev
 inequality for even, $q-$homogeneous $(q>1)$, strictly convex
functions $ C:\mathbb R^n\to [0,\infty)$. In our case,
$C(x)=\frac{\|x\|^{p'}}{p'}.$}]}: {\it One has}
\begin{equation}\label{LS-p-generic-Dolbeau}
    {\bf Ent}_{dx}(|u|^p)=\int_{\mathbb R^n}|u|^p\log |u|^pdx\leq \frac{n}{p}\log\left(\mathcal
L_{p,n}\|\nabla  u\|_{L^p}^p\right),\ \forall u\in W^{1,p}({\mathbb
R^n}),\ \|u\|_{L^p}=1,
\end{equation}
{\it where the best constant \vspace{-0.4cm}
$$
    \mathcal
L_{p,n}={\small
\frac{p}{n}\left(\frac{p-1}{e}\right)^{p-1}\left(\omega_n{\Gamma\left(\frac{n}{p'}+1\right)}\right)^{-\frac{p}{n}}}
$$
is achieved  by the family of functions}
$$
l_p^\lambda(x)={\small
\lambda^\frac{n}{pp'}\omega_n^{-\frac{1}{p}}\Gamma\left(\frac{n}{p'}+1\right)^{-\frac{1}{p}}
e^{-\frac{\lambda}{p}\|x\|^{p'}},\lambda>0};$$

  \item {\bf Limit case II} ($\alpha\to 0$) {\rm [see \cite[p. 320]{CE-N-Villani}]:} {\it
One has
\begin{equation}\label{Faber-Krahn}
   \|u\|_{L^1}\leq \mathcal F_{p,n}\|\nabla  u\|_{L^{p}}|{\rm supp}(u)|^{1-\frac{1}{p^\star}}, \ \forall u\in \dot W^{1,p}({\mathbb
R^n})\ 
\end{equation}
 and the best constant
$$\mathcal F_{p,n}=\lim_{\alpha\to 0} {\mathcal N}_{\alpha,p,n}=n^{-\frac{1}{p}}\omega_n^{-\frac{1}{n}}(p'+n)^{-\frac{1}{p'}}$$
is achieved  by the family of functions
$$f_p^\lambda(x)=\lim_{\alpha\to
0}h_{\alpha,p}^\lambda(x)=(\lambda-\|x\|^{p'})_+,\ x\in \mathbb
R^n,$$ where {\rm supp}$(u)$ stands for the support of $u$ and
$|{\rm supp}(u)|$ is its Lebesgue measure.}
\end{itemize}

\subsection{Statement of main results}
\label{sect1.2} As we already pointed out, the primordial purpose of
the present  paper is to establish fine topological properties of
metric measure spaces curved in the sense of Lott-Sturm-Villani
which support Gagliardo-Nirenberg-type inequalities. In fact, the
metric spaces we are working on are supposed to satisfy the
cur\-vature-dimension condition $\textsf{CD}(K,n)$ for some $K\geq
0$ and $n\geq 2,$ introduced
 by Lott and Villani \cite{Lott-Villani} and Sturm
\cite{Sturm1, Sturm2}; see \S \ref{sect-2} for its formal
definition.

\subsubsection{Volume non-collapsing on metric measure spaces}

%

Let $({M},\textsf{d},\textsf{m})$ be a metric measure space (with a
strictly positive Borel measure $\textsf{m}$) 
and ${\rm Lip}_0(M)$ be the space of Lipschitz
functions with compact support on $M$. For  $u\in {\rm Lip}_0(M)$,
let
\begin{equation}\label{local-constant}
    |\nabla u|_{\textsf{{d}}}(x):=\limsup_{y\to x}
\frac{|u(y)-u(x)|}{\textsf{d}(x,y)}, \ x\in M.
\end{equation}
 Note that $x\mapsto |\nabla u|_\textsf{d}(x)$ is Borel measurable on $M$ for $u\in
{\rm Lip}_0(M)$.

As before, let $n\geq 2$ be an integer, $p\in (1,n)$ and $\alpha\in
(0,\frac{n}{n-p}]\setminus \{1\}$. Throughout this section we assume
that the {\it lower $n-$density of the measure {\rm $\textsf{m}$} at
a point $x_0\in M$ is unitary}, i.e.,
\begin{description}
\item[{\rm $({\bf D})^n_{x_0}$}]
 $\displaystyle\liminf_{\rho\to
 0}\frac{\textsf{m}( B_{x_0}(\rho))}{\omega_n\rho^n}=1$,
\end{description}
where $B_{x}(r)=\{y\in M:\textsf{{d}}(x,y)<r\}.$

Throughout the whole paper, we shall keep  the notations from
Theorems A and B (i.e., the four best constants from the
Gagliardo-Nirenberg inequalities on normed spaces and the numbers
$\theta$ and $\gamma$ from (\ref{theta-best}) and
(\ref{gamma-best}), respectively); the Lebesgue spaces $L^p$ are
defined on the measure space $({M},\textsf{m})$. We now are the
position to state our quantitative, globally non-collapsing volume
growth results:

\begin{theorem}\label{main-theorem-Gagliardo} {\rm [Gagliardo-Nirenberg inequalities]} Let $({M},\textsf{\emph{d}},\textsf{\emph{m}})$ be a proper metric measure
space  which satisfies the curvature-dimension condition
$\textsf{\emph{CD}}(K,n)$ for some $K\geq 0$ and $n\geq 2$. Let
$p\in (1,n)$ and assume that {\rm $({\bf D})^n_{x_0}$} holds for
some $x_0\in M$. Then the following statements hold:
\begin{itemize}
  \item[\rm (i)] If $1<\alpha\leq \frac{n}{n-p}$ and the inequality
$$
    \|u\|_{L^{\alpha p}}\leq \mathcal C
    \||\nabla u|_{\sf{{d}}}\|_{L^p}^{\theta}\|u\|_{L^{\alpha(p-1)+1}}^{1-\theta},\
    \forall u\in {\rm Lip}_0(M) \eqno{{\bf (GN1)}_{\mathcal C}^{\alpha,p}}
$$
holds for some $\mathcal C\geq \mathcal G_{\alpha,p,n}$, then $K=0$
and
$${\sf{m}}(B_{x}(\rho))\geq\left(\frac{\mathcal G_{\alpha,p,n}}{\mathcal C}\right)^\frac{n}{\theta}\omega_n\rho^n\ \  for\ all\ x\in M\ and\
\rho\geq 0.$$
\item[\rm (ii)]
If $0<\alpha<1$ and the inequality
$$
    \|u\|_{L^{\alpha(p-1)+1}}\leq \mathcal C
    \||\nabla u|_{\sf{{d}}}\|_{L^p}^{\gamma}\|u\|_{L^{\alpha p}}^{1-\gamma},\
    \forall u\in {\rm Lip}_0(M)
 \eqno{{\bf (GN2)}_{\mathcal C}^{\alpha,p}}
$$
holds for some $\mathcal C\geq \mathcal N_{\alpha,p,n}$, then $K=0$
and
$${\sf m}(B_{x}(\rho))\geq\left(\frac{\mathcal N_{\alpha,p,n}}{\mathcal C}\right)^\frac{n}{\gamma}\omega_n\rho^n\ \  for\ all\ x\in M\ and\
\rho\geq 0.$$
\end{itemize}
\end{theorem}

In the limit case $\alpha\to 1$, we can state
\begin{theorem}\label{main-theorem-log-Sobolev} {\rm [$L^p-$logarithmic Sobolev inequality]} Under the same assumptions as in Theorem \ref{main-theorem-Gagliardo},
if
$$
   {\bf Ent}_{d{\sf m}}(|u|^p)= \int_{M}|u|^p\log |u|^pd{\sf m}\leq \frac{n}{p}\log\left(\mathcal
C\||\nabla u|_{\sf{{d}}}\|_{L^p}^p \right),\ \forall u\in {\rm
Lip}_0(M),\ \|u\|_{L^p}=1 \eqno{{\bf (LS)}_{\mathcal C}^{p}}
$$
holds for some $\mathcal C\geq \mathcal L_{p,n},$ then $K=0$ and
$$\textsf{\emph{m}}( B_x(\rho))\geq \left(
\frac{\mathcal L_{p,n}}{\mathcal C}
\right)^\frac{n}{p}\omega_n\rho^n\ \  for\ all\ x\in M\ and\
\rho\geq 0.$$
\end{theorem}

\noindent   In the remaining limit case $\alpha\to 0$, one can prove

\begin{theorem}\label{main-theorem-Faber-Krahn} {\rm [Faber-Krahn-type inequality]} Under the same assumptions as in Theorem \ref{main-theorem-Gagliardo}, if
$$
   \|u\|_{L^1}\leq \mathcal C\||\nabla u|_{\sf{{d}}}\|_{L^p}{\sf m}({\rm supp}(u))^{1-\frac{1}{p^\star}},\ \forall u\in {\rm
Lip}_0(M) \eqno{{\bf (FK)}_{\mathcal C}^{p}}
$$
holds for some $\mathcal C\geq \mathcal F_{p,n}$, then $K=0$ and
$$\textsf{\emph{m}}( B_x(\rho))\geq\left(
\frac{\mathcal F_{p,n}}{\mathcal C} \right)^{n}\omega_n\rho^n\ \
for\ all\ x\in M\ and\ \rho\geq 0.$$
\end{theorem}

\noindent Some remarks are in order.

\begin{remark}\rm \label{remark-bonnet}
(a) The proofs of Theorems
\ref{main-theorem-Gagliardo}-\ref{main-theorem-Faber-Krahn} are {\it
synthetic} where
 we shall exploit some basic features of metric
measure spaces satisfying the $\textsf{CD}(K,n)$ condition
 (such as generalized Bonnet-Myers and
Bishop-Gromov comparison inequalities) and direct constructions.
Although the lines of the proofs of these results are similar, our
arguments
 require different technics, deeply depending on the {\it shape}
 of certain test functions whose profiles come from the family of extremals
in normed spaces (cf. Theorems A \& B). Note that instead of the
$\textsf{CD}(K,n)$ condition it is enough to consider the slightly
weaker {\it measure contraction property} $\textsf{MCP}(K,n)$, see
Ohta \cite{Ohta-CMH}.


 (b) The case  $p=2$ and $\alpha=\frac{n}{n-2}$ ($n\geq 3$) is contained in  Krist\'aly and Ohta \cite{Kri-Ohta}, where the authors studied
 Caffarelli-Kohn-Nirenberg inequalities on metric measure spaces. We notice that the roots of Theorem \ref{main-theorem-Gagliardo} (i)
 on Riemannian manifolds with non-negative Ricci curvature can be found
 in  do Carmo and Xia \cite{doCarmo-Xia}, Ledoux \cite{Ledoux-CAG} and Xia \cite{Xia-JFA}.

 (c) The generalized Bishop-Gromov inequality and density assumption $({\bf D})^n_{x_0}$
 imply  ${\sf m}( B_{x_0}(\rho))\leq\omega_n\rho^n\ \  {\rm for\ all}\
\rho\geq 0.$ In particular, the latter
 inequality and the
conclusions of Theorems
\ref{main-theorem-Gagliardo}-\ref{main-theorem-Faber-Krahn} imply
the Ahlfors $n-$regularity at the point $x_0$; therefore, the
Hausdorff dimension of $(M,{\sf d})$ is precisely $n.$

 (d) $({\bf D})^n_{x_0}$ clearly holds for every point $x_0$ on $n-$di\-men\-sional Riemannian and Finsler manifolds endowed with the canonical Busemann-Hausdorff measure.
\end{remark}

\subsubsection{Applications: rigidity results in smooth
settings}\label{sect1.3}
 Having fine volume growth estimates in Theorems
\ref{main-theorem-Gagliardo}-\ref{main-theorem-Faber-Krahn},
important {\it rigidity} results can be deduced in the context of
 Riemannian and Finsler manifolds supporting Gagliardo-Nirenberg-type
 inequalities.

In order to state such results, let $(M,g)$ be an $n-$dimensional
complete Riemannian manifold with non-negative Ricci curvature
$(n\geq 2)$ endowed with its canonical volume form $dv_g$. Let
$\alpha_{MP}(k,n)\in (0,1]$ be the so-called {\it Munn-Perelman
constant} for every $k=1,...,n,$ see Munn \cite{Munn-JGA}.  In fact,
based on the double induction argument of Perelman
\cite{Perelman-JAMS}, Munn determined explicit lower bounds for the
volume growth  in terms of the constant $\alpha_{MP}(k,n)$ which
guarantee the triviality of the $k$-th
homotopy group $\pi_k(M)$ of $(M,g);$ see details in \S \ref{sect-4}. 

 For sake of simplicity,  we restrict here our attention to the
$L^p-$logarithmic Sobolev inequality$({\bf LS})_{\mathcal C}^p$ on
$(M,g)$ by proving that once $\mathcal C>0$ is closer and closer to
the optimal Euclidean constant $\mathcal L_{p,n}$, the manifold
$(M,g)$ approaches topologically more and more to the Euclidean
space $\mathbb R^n.$

\begin{theorem}\label{theorem-Riemann-log-Sobolev}
Let $(M,g)$ be an $n-$dimensional complete Riemannian manifold
 with non-negative Ricci
curvature $(n\geq 2)$ and assume the $L^p-$logarithmic Sobolev
inequality $({\bf LS})_{\mathcal C}^p$ holds on $(M,g)$ for some
$p\in (1,n)$ and $\mathcal C>0$. Then the following assertions hold:
\begin{itemize}
 \item[{\rm (i)}] $\mathcal C\geq \mathcal L_{p,n};$
  \item[{\rm (ii)}] The order of the fundamental group $\pi_1(M)$ is bounded
  above by $\left(
\frac{\mathcal C}{\mathcal L_{p,n}} \right)^\frac{n}{p};$
 \item[{\rm (iii)}] If $\mathcal C < \alpha_{MP}(k_0,n)^{-\frac{p}{n}}\mathcal L_{p,n}$ for some $k_0\in \{1,...,n\}$ then
 $\pi_1(M)=...=\pi_{k_0}(M)=0;$
 \item[{\rm (iv)}] If $\mathcal C < \alpha_{MP}(n,n)^{-\frac{p}{n}}\mathcal L_{p,n}$  then
 $M$ is contractible$;$
  \item[{\rm (v)}]
  $\mathcal C=\mathcal L_{p,n}$   if and only if $(M,g)$ is isometric to the Euclidean space $\mathbb
  R^n.$
\end{itemize}
\end{theorem}

\begin{remark}\rm
(a) Theorem \ref{theorem-Riemann-log-Sobolev} (v) answers an open
question of Xia \cite{Xia} for generic $p\in (1,n)$. For $p=2$ the
latter equivalence is well known by using sharp analytic estimates
for the heat kernel on complete Riemannian manifolds with
non-negative Ricci curvature; see  Bakry, Concordet and Ledoux
\cite{BCL}, Ni \cite{Ni}, and Li \cite{Peter_Li}. Details are
presented in \S \ref{sect-riemann-3} (see Remark \ref{rem-3.1}).

(b) The conclusion $\mathcal C\geq \mathcal L_{p,n}$ in Theorem
\ref{theorem-Riemann-log-Sobolev} (i) is in a perfect concordance
with the assumption of Theorem \ref{main-theorem-log-Sobolev}.
Analogous statements hold for the other  Gagliardo-Nirenberg
inequalities.

(c) Similar results to Theorem \ref{theorem-Riemann-log-Sobolev} can
be stated also for Gagliardo-Nirenberg inequalities ${\bf
(GN1)}_{\mathcal C}$ and ${\bf (GN2)}_{\mathcal C}$, and Faber-Krahn
inequality $({\bf FK})_{\mathcal C}$ with trivial modifications. In
particular, we have:
\end{remark}

\begin{corollary}\label{corollary-Riemann} {\rm [Optimality vs. flatness]} Let $(M,g)$ be an $n(\geq 2)-$dimensional complete Riemannian manifold  with non-negative Ricci
curvature. The following statements are equivalent:
\begin{itemize}
\item[{\rm (i)}]  ${\bf (GN1)}_{\mathcal G_{\alpha,p,n}}^{\alpha, p}$ holds on $(M,g)$ for some $p\in (1,n)$ and $\alpha\in
(1,\frac{n}{n-p}];$
\item[{\rm (ii)}]  ${\bf (GN2)}_{\mathcal N_{\alpha,p,n}}^{\alpha, p}$ holds on $(M,g)$ for some $p\in (1,n)$ and $\alpha\in
(0,1);$
  \item[{\rm (iii)}]  $({\bf
LS})_{\mathcal L_{p,n}}^p$ holds on $(M,g)$ for some $p\in (1,n);$
 \item[{\rm (iv)}]  $({\bf
FK})_{\mathcal F_{p,n}}^p$ holds on $(M,g)$ for some $p\in (1,n);$
  \item[{\rm (v)}] $(M,g)$ is isometric to the Euclidean space
$\mathbb R^n$.
\end{itemize}
\end{corollary}

\begin{remark}\rm  (a) The equivalence (i)$\Leftrightarrow$(v) in Corollary \ref{corollary-Riemann} is
precisely the main result of Xia \cite{Xia-JFA}.

(b) A similar rigidity result to Corollary \ref{corollary-Riemann}
can be stated on reversible Finsler manifolds endowed with the
natural Busemann-Hausdoff measure $dV_F$ of $(M,F)$; roughly
speaking, we can replace the notions 'Riemannian' and 'Euclidean' in
Corollary \ref{corollary-Riemann} by the notions 'Berwald' and
'Minkowski', respectively (see Theorem \ref{theorem-Finsler-1}). The
latter notions will be introduced in \S\ref{sect-finsler}.
\end{remark}

\noindent {\bf Notations.}  When no confusion arises,
 $\|\cdot\|_{L^p}$ abbreviates: (a) $\|\cdot\|_{L^p{(M,d\sf
m})}$ on the metric measure space $(M,{\sf d},{\sf m})$; (b)
$\|\cdot\|_{L^p(M,dv_g)}$ on the Riemannian manifold $(M,g)$ where
$dv_g$ stands for the canonical Riemannian measure on $(M,g)$; (c)
$\|\cdot\|_{L^p(M,dV_F)}$ on the Finsler manifold $(M,F)$ where
$dV_F$ denotes the Busemann-Hausdoff measure on $(M,F)$; and (d)
$\|\cdot\|_{L^p(\mathbb R^n,dx)}$ on the Euclidean/normed space
$\mathbb R^n$ where $dx$ is the usual Lebesgue measure,
respectively. When $A$ is not the whole space we are working on, we
shall use the notation
 $\|u\|_{L^p(A)}$ for the $L^p-$norm of the function $u:A\to \mathbb R$.


%
%


\section{\sc Volume non-collapsing via Gagliardo-Nirenberg inequalities} \label{sect-2}

Before the presentation of the proofs of Theorems
\ref{main-theorem-Gagliardo}-\ref{main-theorem-Faber-Krahn}, we
recall for completeness some notions and results from Lott and
Villani \cite{Lott-Villani} and Sturm \cite{Sturm1, Sturm2}, which
are indispensable in our arguments.

Let $({M},\textsf{d},\textsf{m})$ be a metric measure space, i.e.,
$({M},\textsf{d})$ is a complete separable metric space and
$\textsf{m}$ is a locally finite measure on $M$ endowed with its
Borel $\sigma-$algebra. In the sequel, we assume that the measure
$\textsf{m}$ on $M$ is strictly positive, i.e.,
supp$[\textsf{m}]=M.$ As usual, $\mathcal P_2(M,\textsf{d})$ is the
$L^2-$Wasserstein space of probability measures on $M$, while
$\mathcal P_2(M,\textsf{d},\textsf{m})$ will denote the subspace of
$\textsf{m}-$absolutely continuous measures.
$({M},\textsf{d},\textsf{m})$ is said to be proper if every bounded
and closed subset of $M$ is compact.

For a given number  $N\geq 1,$ the {\it R\'enyi entropy functional}
$S_N(\cdot|\textsf{m}):\mathcal P_2(M,\textsf{d})\to \mathbb R$ with
respect to the measure $\textsf{m}$ is defined by
$S_N(\mu|\textsf{m})=-\int_M \rho^{-\frac{1}{N}}d\mu,$
 $\rho$ being the density of $\mu^c$ in
$\mu=\mu^c+\mu^s=\rho \textsf{m}+\mu^s$, where $\mu^c$ and $\mu^s$
represent the absolutely continuous and singular parts of $\mu\in
\mathcal P_2(M,\textsf{d}),$ respectively.

Let $K,N\in \mathbb R$ be two numbers with $K\geq 0$ and $N\geq 1$.
For every $t\in [0,1]$ and $s\geq 0$, let {\small
$$\tau_{K,N}^{(t)}(s)=\left\{
\begin{array}{lll}
+\infty, & {\rm if} & Ks^2\geq (N-1)\pi^2;
\\ t^\frac{1}{N}\left(\sin\left(\sqrt{\frac{K}{N-1}}ts\right)\big/\sin\left(\sqrt{\frac{K}{N-1}}s\right)\right)^{1-\frac{1}{N}},& {\rm if} &
0<Ks^2<(N-1)\pi^2;\\
t, & {\rm if} & Ks^2=0.
\end{array}
\right.$$}

 We say that $({M},\textsf{d},\textsf{m})$ satisfies the {\it
curvature-dimension condition} $\textsf{CD}(K,N)$ if for each
$\mu_0,\mu_1\in \mathcal P_2(M,\textsf{d},\textsf{m})$ there exists
an optimal coupling $\gamma$ of $\mu_0,\mu_1$ and a geodesic
$\Gamma:[0,1]\to \mathcal P_2(M,\textsf{d},\textsf{m})$ joining
$\mu_0$ and $\mu_1$ such that
$$S_{N'}(\Gamma(t)|\textsf{m})\leq - \int_{M\times M}[\tau_{K,N'}^{(1-t)}(\textsf{d}(x_0,x_1))\rho_0^{-\frac{1}{N'}}(x_0)+
\tau_{K,N'}^{(t)}(\textsf{d}(x_0,x_1))\rho_1^{-\frac{1}{N'}}(x_1)]d\gamma(x_0,x_1)$$
for every $t\in [0,1]$ and $N'\geq N$, where $\rho_0$ and $\rho_1$
are the densities of $\mu_0$ and $\mu_1$ with respect to
$\textsf{m}$. Clearly, when $K=0$, the above inequality reduces to
the the geodesic convexity of $ S_{N'}(\cdot|\textsf{m})$ on the
$L^2-$Wasserstein space $\mathcal P_2(M,\textsf{d},\textsf{m})$.

It is well known that $\textsf{CD}(K,n)$ holds on a complete
Riemannian manifold $(M,g)$ endowed with the Riemannian volume
element $dv_g$ if and only if its Ricci curvature $\geq K$ and
dim$(M)\leq n.$

Let $B_{x}(r)=\{y\in M:\textsf{{d}}(x,y)<r\}$. In the sequel we
shall exploit properties which are resumed in the following results.

\begin{theorem} {\rm (see \cite{Sturm2})}\label{CD-tetel}
Let $({M},\textsf{\emph{d}},\textsf{\emph{m}})$ be a metric measure
space with strictly positive measure $\textsf{\emph{m}}$ satisfying
the curvature-dimension condition $\textsf{\emph{CD}}(K,N)$ for some
$K\geq 0$ and $N>1$. Then every bounded set $S\subset M$ has finite
$\textsf{\emph{m}}-$measure and the metric spheres $\partial
B_{x}(r)$ have zero $\textsf{\emph{m}}-$measures. Moreover, one has:
\begin{itemize}
  \item[{\rm (i)}] {\rm [Generalized Bonnet-Myers theorem]} If $K>0,$
  then
  $M=${\rm supp}$[{\sf m}]$ is compact and has diameter less than or equal to $\sqrt{\frac{N-1}{K}}\pi.$
  \item[{\rm (ii)}] {\rm [Generalized Bishop-Gromov inequality]} If
  $K=0,$ then for every $R>r>0$ and $x\in M,$
  $$\frac{\textsf{\emph{m}}( B_x(r))}{r^N}\geq\frac{\textsf{\emph{m}}( B_x(R))}{R^N}.$$
\end{itemize}
\end{theorem}

\begin{lemma}\label{lemma-independence-of-point} Let $({M},\textsf{\emph{d}},\textsf{\emph{m}})$ be a metric measure
space which satisfies the curvature-dimension condition
$\textsf{\emph{CD}}(0,n)$ for some $n\geq 2$. If
\begin{equation}\label{ell-ell}
    \ell_\infty^{x_0}:=\limsup_{\rho\to \infty}\frac{{\sf m}(
    B_{x_0}(\rho))}{\omega_n\rho^n}\geq a
\end{equation}
    for some $x_0\in M$ and $a>0$, then
    $${\sf m}(
    B_{x}(\rho))\geq a \omega_n\rho^n,\ \forall x\in M,\ \rho\geq 0.$$
\end{lemma}

{\it Proof.} Let us fix $x\in M$ and $\rho>0$; then we have
\begin{eqnarray*}
   \frac{\textsf{m}(
    B_{x}(\rho))}{\omega_n\rho^n} &\geq& \limsup_{r\to \infty}\frac{\textsf{m}(
    B_{x}(r))}{\omega_nr^n}\ \ \ \ \ \ \ \ \ \ \ \ \  \ \ \ \ \ \ \ \ \ \ \  \ \ \ \ \ \ \ \ \ {\rm [Bishop-Gromov\ inequality]} \\
   &\geq&   \limsup_{r\to \infty}\frac{\textsf{m}(
    B_{x_0}(r-\textsf{d}(x_0,x)))}{\omega_nr^n}\ \ \ \ \ \ \ \ \ \ \ \ \  \ \ \ \ \ {\rm [}B_{x}(r)\supset B_{x_0}(r-\textsf{d}(x_0,x)){\rm ]} \\
   &=& \limsup_{r\to \infty}\left(\frac{\textsf{m}(
    B_{x_0}(r-\textsf{d}(x_0,x)))}{\omega_n(r-\textsf{d}(x_0,x))^n}\cdot\frac{(r-\textsf{d}(x_0,x))^n}{r^n}\right)\\&=&
\ell_\infty^{x_0}
\\&\geq&a,\ \ \ \ \ \ \ \
\ \ \ \ \ \ \ \  \ \ \ \ \ \ \ \ \ \ \ \ \ \ \ \ \ \ \ \  \  \ \ \ \
\ \ \ \ \ \ \ \ \ \ \ \  {\rm [cf.\
 (\ref{ell-ell})]}
\end{eqnarray*}
which concludes the proof. \hfill $\square$\\

\noindent We are now in the position to prove our volume
non-collapsing results.

\subsection{Cases $\alpha>1$ \& $0<\alpha<1$: usual Gagliardo-Nirenberg inequalities
} In this subsection we present the proof of Theorem
\ref{main-theorem-Gagliardo} by distinguishing two cases:\\

 {\it Proof of Theorem {\rm \ref{main-theorem-Gagliardo} (i)}:
the case $1<\alpha\leq \frac{n}{n-p}$.} In this part we follow the
line of \cite{Kri-Ohta}; the proof is divided into several steps. We
clearly may assume that $\mathcal C>\mathcal G_{\alpha,p,n}$ in
${\bf (GN1)}_{\mathcal C}^{\alpha,p}$; indeed, if $\mathcal
C=\mathcal G_{\alpha,p,n}$ we can consider the subsequent arguments
for $\mathcal C:=\mathcal G_{\alpha,p,n}+\varepsilon$ with small
$\varepsilon>0$ and then take $\varepsilon\to 0^+.$

{\sc Step 1 $(K=0)$.} If we assume that $K>0$ then the generalized
Bonnet-Myers theorem (see Theorem \ref{CD-tetel} (i)) implies that
$M$ is compact and $\textsf{m}(M)$ is finite. Taking the constant
map $u(x)={\textsf{m}(M)}$ in ${\bf (GN1)}_{\mathcal C}^{\alpha,p}$
as a test function, one gets a contradiction. Therefore, $K=0.$

{\sc Step 2 ({\it ODE from the optimal Euclidean Gagliardo-Nirenberg
inequality I}).} We consider the optimal Gagliardo-Nirenberg
inequality (\ref{Villani-1}) in the particular case when the norm is
precisely the Euclidean norm $|\cdot|.$ After a simple rescaling,
one can see that the function $x\mapsto
(\lambda+|x|^{p'})^\frac{1}{1-\alpha},$ $\lambda>0,$ is a family of
extremals in (\ref{Villani-1}); therefore, we have the following
first order
 ODE
\begin{equation}\label{ODE-1}
   \left(\frac{1-\alpha}{\alpha(p-1)+1}h_G'(\lambda)\right)^\frac{1}{\alpha
   p}=\mathcal G_{\alpha,p,n}\left(\frac{p'}{\alpha-1}\right)^\theta
   \left(h_G(\lambda)+\frac{\alpha-1}{\alpha(p-1)+1}\lambda
   h_G'(\lambda)\right)^\frac{\theta}{p}h_G(\lambda)^\frac{1-\theta}{\alpha(p-1)+1},
\end{equation}
where $h_G:(0,\infty)\to \mathbb R$ is given by
$$ h_G(\lambda) = \int_{\mathbb
R^n}\left(\lambda+|x|^{p'}\right)^\frac{\alpha(p-1)+1}{1-\alpha}dx,\
\lambda>0.$$ For further use, we shall represent the function $h_G$
in two different ways, namely \begin{eqnarray}\label{h-G-repres}
 \nonumber h_G(\lambda)
   &=& \omega_n\frac{n}{p'}{\sf
B}\left(\frac{\alpha(p-1)+1}{\alpha-1}-\frac{n}{p'},\frac{n}{p'}\right)
\lambda^{\frac{\alpha(p-1)+1}{1-\alpha}+\frac{n}{p'}}\\
&=&\int_0^\infty \omega_n\rho^nf_G(\lambda,\rho)d\rho,
\end{eqnarray}
where
\begin{equation}\label{f-representation}
f_G(\lambda,\rho)=p'\frac{\alpha(p-1)+1}{\alpha-1}\left(\lambda+\rho^{p'}\right)^\frac{\alpha
p}{1-\alpha}\rho^{p'-1}.
\end{equation}

{\it {\sc Step 3} $($Differential inequality from $({\bf
GN1})_{\mathcal C}^{\alpha,p}$$)$.} By the generalized Bishop-Gromov
inequality (see Theorem \ref{CD-tetel} (ii)) and hypothesis {\rm
$({\bf D})^n_{x_0}$} one has that
\begin{equation}\label{x_0koruli_novekedes}
   \frac{\textsf{m}(B_{x_0}(\rho))}{\omega_n\rho^n}\leq \displaystyle\liminf_{r\to
 0}\frac{\textsf{m}(B_{x_0}(r))}{\omega_nr^n}=1,\ \rho>0.
\end{equation}
Inspired by the form of $h_G$, we consider  the function
$w_G:(0,\infty)\to \mathbb R$ defined by
$$ w_G(\lambda) = \int_{M}\left(\lambda+{\sf d}(x_0,x)^{p'}\right)^\frac{\alpha(p-1)+1}{1-\alpha}d{\sf m}(x),\
\lambda>0.$$ By using the layer cake representation, it follows that
$w_G$ is well-defined and of class $C^1$; indeed,
\begin{eqnarray*}
  w_G(\lambda) &=& \int_0^\infty {\sf m}\left(\left\{x\in M:\left(\lambda+{\sf d}(x_0,x)^{p'}\right)^\frac{\alpha(p-1)+1}{1-\alpha}>t\right\}\right)dt \\
   &=& \int_0^\infty  {\sf m}(B_{x_0}(\rho)) f_G(\lambda,\rho)d\rho\ \ \ \ \ \ \ \ \ \ \ \ \  \ \ \ {\rm [change}\ t=\left(\lambda+\rho^{p'}\right)^\frac{\alpha(p-1)+1}{1-\alpha}\ {\rm and\ see\ (\ref{f-representation})}]\\
   &\leq & \int_0^\infty \omega_n\rho^nf_G(\lambda,\rho)d\rho \ \ \ \ \
   \ \ \ \
    \ \ \ \ \ \ \ \ \ \ \ \  \ \ {\rm [see\
    (\ref{x_0koruli_novekedes})]} \\
    &=& h_G(\lambda),
\end{eqnarray*}
thus
\begin{equation}\label{kicsi-hasonlitas}
    0< w_G(\lambda)\leq h_G(\lambda)<\infty,\ \lambda>0.
\end{equation}

For every $\lambda>0$ and $k\in \mathbb N$, we consider the function
$u_{\lambda,k}:M\to \mathbb R$ defined by
\[ u_{\lambda,k}(x)
 =(\min\{0,k-{\sf d}(x_0,x)\}+1)_+
 \left(\lambda+\max\left\{ {\sf d}(x_0,x),k^{-1}\right\}^{p'}\right)^{^\frac{1}{1-\alpha}}. \]
Note that since $({M},\textsf{d},\textsf{m})$ is proper, the set
${\rm supp}(u_{\lambda,k})=\overline {B_{x_0}(k+1)}$ is compact.
Consequently, $u_{\lambda,k}\in {\rm Lip}_0(M)$ for every
$\lambda>0$ and $k\in \mathbb N$; thus we can apply these functions
in $({\bf GN1})_{\mathcal C}^{\alpha,p}$, i.e.,
$$\|u_{\lambda,k}\|_{L^{\alpha p}}\leq \mathcal C
    \||\nabla u_{\lambda,k}|_{\sf{{d}}}\|_{L^p}^{\theta}\|u_{\lambda,k}\|_{L^{\alpha(p-1)+1}}^{1-\theta}.$$
 Moreover,
\[  \lim_{k\to \infty}u_{\lambda,k}(x)= \left(\lambda+{\sf d}(x_0,x)^{p'}\right)^\frac{1}{1-\alpha}=:u_\lambda(x)
. \]  By using the dominated convergence theorem, it turns out from
the above inequality that $u_{\lambda}$ also verifies $({\bf
GN1})_{\mathcal C}^{\alpha,p}$, i.e.,
\begin{equation}\label{u-lambda-egyenlotlenseg}
    \|u_{\lambda}\|_{L^{\alpha p}}\leq \mathcal C
    \||\nabla u_{\lambda}|_{\sf{{d}}}\|_{L^p}^{\theta}\|u_{\lambda}\|_{L^{\alpha(p-1)+1}}^{1-\theta}.
\end{equation}
The non-smooth chain rule gives that
\begin{equation}\label{chain-rule-lipschitz}
|\nabla
u_{\lambda}|_{\sf{{d}}}(x)=\frac{p'}{\alpha-1}\left(\lambda+{\sf
d}(x_0,x)^{p'}\right)^\frac{\alpha }{1-\alpha}{\sf
d}(x_0,x)^{p'-1}|\nabla {\sf d}(x_0,\cdot)|_{\sf{{d}}}(x),\ x\in
M.\end{equation}
 Since ${\sf
d}(x_0,\cdot)$ is $1$-Lipschitz (therefore, $|\nabla {\sf
d}(x_0,\cdot)|_{\sf{{d}}}(x)\leq 1$ for all $x\in M$), due to
(\ref{u-lambda-egyenlotlenseg}), (\ref{chain-rule-lipschitz}) and
the form of the function $w_G$, we obtain the differential
inequality
\begin{equation}\label{diff-egyenlotlenseg-1}
   \left(\frac{1-\alpha}{\alpha(p-1)+1}w_G'(\lambda)\right)^\frac{1}{\alpha
   p}\leq \mathcal C\left(\frac{p'}{\alpha-1}\right)^\theta
   \left(w_G(\lambda)+\frac{\alpha-1}{\alpha(p-1)+1}\lambda
   w_G'(\lambda)\right)^\frac{\theta}{p}w_G(\lambda)^\frac{1-\theta}{\alpha(p-1)+1}.
\end{equation}

{\it {\sc Step 4} $($Comparison of $w_G$ and $h_G$ near the
origin$).$} We claim that
\begin{equation}\label{h_G-w_G-comparison}
    \lim_{\lambda\to 0^+}\frac{w_G(\lambda)}{h_G(\lambda)}=1.
\end{equation}
By hypothesis {\rm $({\bf D})^n_{x_0}$},  for every $\varepsilon>0$
there exists $\rho_\varepsilon>0$ such that
\begin{equation}\label{local-becsles-gomb}
    \textsf{m}(B_{x_0}(\rho))\geq (1-\varepsilon)\omega_n\rho^n \ {\rm
for\ all}\ \rho\in [0,\rho_\varepsilon].
\end{equation}
By (\ref{local-becsles-gomb}), one has that
\begin{eqnarray*}
  w_G(\lambda) &=& \int_0^\infty  {\sf m}(B_{x_0}(\rho)) f_G(\lambda,\rho)d\rho \\
   &\geq&  (1-\varepsilon) \int_0^{\rho_\varepsilon}  \omega_n \rho^n f_G(\lambda,\rho)d\rho=
   (1-\varepsilon)\lambda^{\frac{\alpha(p-1)+1}{1-\alpha}+\frac{n}{p'}} \int_0^{\rho_\varepsilon\lambda^{-\frac{1}{p'}}}  \omega_n \rho^n
   f_G(1,\rho)d\rho.
\end{eqnarray*}
Thus, by the representation (\ref{h-G-repres}) of  $h_G$ and a
change of variables, it turns out that
$$\liminf_{\lambda\to
0^+}\frac{w_G(\lambda)}{h_G(\lambda)}\geq
(1-\varepsilon)\liminf_{\lambda\to
0^+}\frac{\displaystyle\int_0^{\rho_\varepsilon\lambda^{-\frac{1}{p'}}}
\omega_n \rho^n f_G(1,\rho)d\rho}{\displaystyle\int_0^\infty
\omega_n \rho^n f_G(1,\rho)d\rho}= 1-\varepsilon.$$ The above
inequality (with $\varepsilon>0$ arbitrary small) combined with
(\ref{kicsi-hasonlitas}) proves the claim
(\ref{h_G-w_G-comparison}).

{\it {\sc Step 5} $($Global comparison of $w_G$ and $h_G$$).$} We
now claim that
\begin{equation}\label{g-h-global-GN1}
    w_G(\lambda)\geq \left(\frac{\mathcal G_{\alpha,p,n}}{\mathcal
C}\right)^\frac{n}{\theta}h_G(\lambda)=\tilde
h_G(\lambda),\ \lambda>0.
\end{equation}
Since we assumed that $\mathcal C>\mathcal G_{\alpha, p,n},$ by
(\ref{h_G-w_G-comparison}) one has   $$\lim_{\lambda\to
0^+}\frac{w_G(\lambda)}{\tilde h_G(\lambda)}= \left(\frac{\mathcal
C}{\mathcal G_{\alpha,p,n}}\right)^\frac{n}{\theta}>1.$$ Therefore,
there exists $\lambda_0>0$ such that for every $\lambda\in
(0,\lambda_0)$, one has $w_G(\lambda)> \tilde h_G(\lambda).$

By contradiction to (\ref{g-h-global-GN1}), we assume that there
exists $\lambda^\#>0$ such that $w_G(\lambda^\#)< \tilde
h_G(\lambda^\#).$ If
$\lambda^*=\sup\{0<\lambda<\lambda^\#:w_G(\lambda)= \tilde
h_G(\lambda)\}$, then  $0<\lambda_0\leq \lambda^*<\lambda^\#.$ In
particular,
$$w_G(\lambda)\leq \tilde
h_G(\lambda),\ \forall \lambda\in [\lambda^*,\lambda^\#].$$ The
latter relation and the differential inequality
(\ref{diff-egyenlotlenseg-1}) imply that for every $\lambda\in
[\lambda^*,\lambda^\#]$,
 \begin{equation}\label{1111}  \left(\frac{1-\alpha}{\alpha(p-1)+1}w_G'(\lambda)\right)^\frac{1}{\alpha
   \theta}\leq \mathcal
   C^\frac{p}{\theta}\left(\frac{p'}{\alpha-1}\right)^p
   \left(\tilde h_G(\lambda)+\frac{\alpha-1}{\alpha(p-1)+1}\lambda
   w_G'(\lambda)\right)\tilde h_G(\lambda)^\frac{(1-\theta)p}{\theta(\alpha(p-1)+1)}.
\end{equation}
Moreover, since $\tilde h_G(\lambda)=\left(\frac{\mathcal
G_{\alpha,p,b}}{\mathcal C}\right)^\frac{n}{\theta}h_G(\lambda)$,
the ODE in (\ref{ODE-1}) can be equivalently transformed for every
$\lambda>0$ into the equation
\begin{equation} \label{2222} \left(\frac{1-\alpha}{\alpha(p-1)+1}\tilde h_G'(\lambda)\right)^\frac{1}{\alpha
   \theta}= \mathcal
   C^\frac{p}{\theta}\left(\frac{p'}{\alpha-1}\right)^p
   \left(\tilde h_G(\lambda)+\frac{\alpha-1}{\alpha(p-1)+1}\lambda
   \tilde h_G'(\lambda)\right)\tilde
   h_G(\lambda)^\frac{(1-\theta)p}{\theta(\alpha(p-1)+1)}.
\end{equation}
 For $\lambda>0$ fixed we introduce the increasing
function $j_G^\lambda:(0,\infty)\to \mathbb R$ defined by
$$j_G^\lambda(t)=\left(\frac{\alpha-1}{\alpha(p-1)+1}t\right)^\frac{1}{\alpha
   \theta}+\mathcal C^\frac{p}{\theta}\left(\frac{p'}{\alpha-1}\right)^p\frac{\alpha-1}{\alpha(p-1)+1}\lambda\tilde h_G(\lambda)^\frac{(1-\theta)p}{\theta(\alpha(p-1)+1)} t.$$
Relations (\ref{1111}) and (\ref{2222}) can be rewritten into
$$j_G^\lambda(-w_G'(\lambda))\leq \mathcal C^\frac{p}{\theta}\left(\frac{p'}{\alpha-1}\right)^p\tilde h_G(\lambda)^{1+\frac{(1-\theta)p}{\theta(\alpha(p-1)+1)}}=j_G^\lambda(-\tilde h_G'(\lambda)),\ \forall \lambda\in [\lambda^*,\lambda^\#],$$
which implies that $$-w_G'(\lambda)\leq-\tilde h_G'(\lambda),\
\forall \lambda\in [\lambda^*,\lambda^\#],$$ i.e., the function
$\tilde h_G-w_G$ is non-increasing in $[\lambda^*,\lambda^\#]$. In
particular, $0<(\tilde h_G-w_G)(\lambda^\#)\leq (\tilde
h_G-w_G)(\lambda^*)=0$, a contradiction. This concludes the proof of
(\ref{g-h-global-GN1}).

 {\it {\sc Step 6} $($Asymptotic volume growth estimate w.r.t.
 $x_0$$)$.} We claim
that
\begin{equation}\label{asymptotic-GN1}
    \ell_\infty^{x_0}:=\limsup_{\rho\to \infty}\frac{\textsf{m}(
    B_{x_0}(\rho))}{\omega_n\rho^n}\geq \left(\frac{\mathcal
G_{\alpha,p,n}}{\mathcal C}\right)^\frac{n}{\theta}.
\end{equation}
By assuming the contrary, there exists $\varepsilon_0>0$ such that
for some $\rho_0>0$,
$$\frac{\textsf{m}(
    B_{x_0}(\rho))}{\omega_n\rho^n}\leq \left(\frac{\mathcal
G_{\alpha,p,n}}{\mathcal C}\right)^\frac{n}{\theta}-\varepsilon_0,\
\forall \rho\geq \rho_0.$$ By (\ref{g-h-global-GN1}) and from the
latter relation,  we have for every $\lambda>0$ that
\begin{eqnarray*}
  0 &\leq& w_G(\lambda)- \left(\frac{\mathcal G_{\alpha,p,n}}{\mathcal
C}\right)^\frac{n}{\theta}h_G(\lambda) \\
   &=& \int_0^\infty  \left(\frac{{\sf m}(B_{x_0}(\rho))}{\omega_n\rho^n}  -  \left(\frac{\mathcal G_{\alpha,p,n}}{\mathcal
C}\right)^\frac{n}{\theta} \right)\omega_n\rho^n f_G(\lambda,\rho)d\rho\\
   &\leq&\left(1+\varepsilon_0- \left(\frac{\mathcal G_{\alpha,p,n}}{\mathcal
C}\right)^\frac{n}{\theta} \right)\int_0^{\rho_0}
\omega_n\rho^nf_G(\lambda,\rho)d\rho-\varepsilon_0 \int_0^\infty
\omega_n\rho^nf_G(\lambda,\rho)d\rho
\end{eqnarray*}
By using (\ref{h-G-repres}), a suitable rearrangement of the terms
in the above relation shows that
$$\varepsilon_0 \frac{n}{p'}{\sf
B}\left(\frac{\alpha(p-1)+1}{\alpha-1}-\frac{n}{p'},\frac{n}{p'}\right)
\lambda^{1+\frac{n}{p'}}\leq \frac{p'}{n+p'}\left(1+\varepsilon_0-
\left(\frac{\mathcal G_{\alpha,p,n}}{\mathcal
C}\right)^\frac{n}{\theta}
\right)\frac{\alpha(p-1)+1}{\alpha-1}\rho_0^{n+p'}.$$ If we take the
limit $\lambda\to +\infty$ in the last estimate, we obtain a
contradiction. Thus, the claim (\ref{asymptotic-GN1}) is proved and
it remains to apply Lemma
\ref{lemma-independence-of-point}, which concludes the proof of Theorem \ref{main-theorem-Gagliardo} (i).  \hfill $\square$\\

 {\it Proof of Theorem {\rm \ref{main-theorem-Gagliardo} (ii)}:
the case $0<\alpha<1$.} We shall invoke some of the arguments from
the proof of Theorem \ref{main-theorem-Gagliardo} (i), emphasizing
that subtle differences arise due to the 'dual' nature of the
Gagliardo-Nirenberg inequalities ${\bf (GN1)}_{\mathcal
C}^{\alpha,p}$ and ${\bf (GN2)}_{\mathcal C}^{\alpha,p}$,
respectively. As before, we may assume that the inequality ${\bf
(GN2)}_{\mathcal C}^{\alpha,p}$ holds with $\mathcal C>\mathcal
N_{\alpha,p,n}$.

{\sc Step 1.} The fact that $K=0$ works similarly as in Theorem
\ref{main-theorem-Gagliardo} (i).

{\sc Step 2.}  Since $x\mapsto
\left(\lambda^{p'}-|x|^{p'}\right)_+^\frac{1}{1-\alpha}$ is an
extremal function in (\ref{Villani-2}) for every $\lambda>0$, we
obtain  the ODE
\begin{eqnarray}\label{1111-GN2}
\nonumber
 h_N(\lambda)^\frac{1}{\alpha(p-1)+1}  &=& \mathcal N_{\alpha,p,n}\left(\frac{p'}{1-\alpha}\right)^\gamma\left(-h_N(\lambda)+\frac{1-\alpha}{p'(\alpha(p-1)+1)}\lambda h_N'(\lambda)\right)^\frac{\gamma}{p}\times \\
   && \times \left(\frac{1-\alpha}{p'(\alpha(p-1)+1)}
\lambda^{1-p'}h_N'(\lambda)\right)^\frac{1-\gamma}{\alpha p},
\end{eqnarray}
where  the function $h_N:(0,\infty)\to \mathbb R$ is defined by
$$h_N(\lambda)=\int_{\mathbb R^n}\left(\lambda^{p'}-|x|^{p'}\right)_+^\frac{\alpha(p-1)+1}{1-\alpha}dx,\
\lambda>0.$$ It is clear that $h_N$  is well-defined, of class $C^1$
and can be represented as
$$h_N(\lambda)=\omega_n\frac{n}{p'}{\sf
B}\left(\frac{\alpha(p-1)+1}{1-\alpha}+1,\frac{n}{p'}\right)
\lambda^{\frac{\alpha p p'}{1-\alpha}+n+p'}=\int_0^\lambda \omega_n
\rho^n f_N(\lambda,\rho)d\rho,$$ where
\begin{equation}\label{f-representation-GN2}
f_N(\lambda,\rho)=p'\frac{\alpha(p-1)+1}{1-\alpha}\left(\lambda^{p'}-\rho^{p'}\right)^\frac{\alpha
p}{1-\alpha}\rho^{p'-1},\  {\rm for\ every}\ \lambda>0\ {\rm and}\
\rho\in (0,\lambda).
\end{equation}

{\sc Step 3.} Let $w_N:(0,\infty)\to \mathbb R$ be the function
defined by
$$w_N(\lambda)=\int_{M}\left(\lambda^{p'}-{\sf d}(x_0,x)^{p'}\right)_+^\frac{\alpha(p-1)+1}{1-\alpha}d{\sf m}(x),\
\lambda>0,$$ where $x_0\in M$ is from {\rm $({\bf D})^n_{x_0}$}. By
the layer cake representation and relations
(\ref{x_0koruli_novekedes}) and (\ref{f-representation-GN2}), $w_N$
is well-defined, positive, of class $C^1$ and
\begin{equation}\label{kicsi-hasonlitas-GN2}
    0< w_N(\lambda)=\int_0^\lambda  {\sf m}(B_{x_0}(\rho)) f_N(\lambda,\rho)d\rho\leq \int_0^\lambda  \omega_n\rho^n f_N(\lambda,\rho)d\rho =h_N(\lambda)<\infty,\ \lambda>0.
\end{equation}
Since  $u_\lambda= \left(\lambda^{p'}-{\sf
d}(x_0,\cdot)^{p'}\right)_+^\frac{1}{1-\alpha}$ is a Lipschitz
function on $M$ with compact support $\overline{B_{x_0}(\lambda)}$,
it belongs to Lip$_0(M)$. Therefore, we may apply $u_\lambda$ in
${\bf (GN2)}_{\mathcal C}^{\alpha,p}$; a similar reasoning as in
(\ref{chain-rule-lipschitz}) leads to  the differential inequality
\begin{eqnarray}\label{2222-GN2}
\nonumber
 w_N(\lambda)^\frac{1}{\alpha(p-1)+1}  &\leq & \mathcal C\left(\frac{p'}{1-\alpha}\right)^\gamma\left(-w_N(\lambda)+\frac{1-\alpha}{p'(\alpha(p-1)+1)}\lambda w_N'(\lambda)\right)^\frac{\gamma}{p}\times \\
   && \times \left(\frac{1-\alpha}{p'(\alpha(p-1)+1)}
\lambda^{1-p'}w_N'(\lambda)\right)^\frac{1-\gamma}{\alpha p},\
\lambda>0.
\end{eqnarray}

{\sc Step 4.} For an arbitrarily fixed $\varepsilon>0$, let
$\rho_\varepsilon>0$ from (\ref{local-becsles-gomb}). If
$0<\lambda<\rho_\varepsilon$, one has that
$$w_N(\lambda)=\int_0^\lambda  {\sf m}(B_{x_0}(\rho)) f_N(\lambda,\rho)d\rho\geq (1-\varepsilon)\int_0^\lambda  \omega_n\rho^n f_N(\lambda,\rho)d\rho =(1-\varepsilon)h_N(\lambda).$$
Consequently, the latter relation together with
(\ref{kicsi-hasonlitas-GN2}) implies that
\begin{equation}\label{h_G-w_G-comparison-GN2}
    \lim_{\lambda\to 0^+}\frac{w_N(\lambda)}{h_N(\lambda)}=1.
\end{equation}

{\sc Step 5.} We shall prove that
\begin{equation}\label{g-h-global-GN2}
    w_N(\lambda)\geq \left(\frac{\mathcal N_{\alpha,p,n}}{\mathcal
C}\right)^\frac{n}{\gamma}h_N(\lambda)=\tilde
h_N(\lambda),\ \lambda>0.
\end{equation}
By (\ref{h_G-w_G-comparison-GN2}) one has
$$\lim_{\lambda\to 0^+}\frac{w_N(\lambda)}{\tilde h_N(\lambda)}=
\left(\frac{\mathcal C}{\mathcal
N_{\alpha,p,n}}\right)^\frac{n}{\gamma}>1,$$ which implies the
existence of a number $\lambda_0>0$ such that $w_N(\lambda)> \tilde
h_N(\lambda)$ for every $\lambda\in (0,\lambda_0)$.

We assume by contradiction that there exists $\lambda^\#>0$ such
that $w_N(\lambda^\#)< \tilde h_N(\lambda^\#).$ If
$\lambda^*=\sup\{0<\lambda<\lambda^\#:w_N(\lambda)= \tilde
h_N(\lambda)\}$, then  $0<\lambda_0\leq \lambda^*<\lambda^\#$ and
\begin{equation}\label{ellent}
w_N(\lambda)\leq \tilde h_N(\lambda),\ \forall \lambda\in
[\lambda^*,\lambda^\#].
\end{equation}
 For every $\lambda>0$, let
$j_N^\lambda:\left(\frac{p'(\alpha(p-1)+1)}{(1-\alpha)\lambda},\infty\right)\to
\mathbb R$ be the function defined by
$$j_N^\lambda(t)=\mathcal C\left(\frac{p'}{1-\alpha}\right)^\gamma\left(-1+\frac{1-\alpha}{p'(\alpha(p-1)+1)}\lambda t\right)^\frac{\gamma}{p}\left(\frac{1-\alpha}{p'(\alpha(p-1)+1)}
\lambda^{1-p'}t\right)^\frac{1-\gamma}{\alpha p}.$$ It is clear that
$j_N^\lambda$ is well-defined, positive and increasing. A direct
computation yields that both values $(\log
w_N)'(\lambda)=\frac{w'_N(\lambda)}{w_N(\lambda)}$ and $(\log \tilde
h_N)'(\lambda)=\frac{\tilde h'_N(\lambda)}{\tilde h_N(\lambda)}$ are
greater than $\frac{p'(\alpha(p-1)+1)}{(1-\alpha)\lambda}$ for every
$\lambda>0$. Taking into account (\ref{gamma-best}), we have
$$\frac{1}{\alpha(p-1)+1}-\frac{\gamma}{p}-\frac{1-\gamma}{\alpha
p}=-\frac{\gamma}{n};$$ therefore,  if we divide the inequality
(\ref{2222-GN2}) by
$w_N(\lambda)^{\frac{\gamma}{p}+\frac{1-\gamma}{\alpha p}}$, we
obtain that
\begin{equation}\label{w-elsosegben}
w_N(\lambda)^{-\frac{\gamma}{n}}\leq j_N^\lambda\left((\log
w_N)'(\lambda)\right),\ \forall \lambda>0.
\end{equation}
 In
a similar manner, by  $\tilde h_N(\lambda)=\left(\frac{\mathcal
N_{\alpha,p,n}}{\mathcal C}\right)^\frac{n}{\gamma}h_N(\lambda)$ and
relation (\ref{1111-GN2}), we have that
\begin{equation}\label{h-elsosegben}
\tilde h_N(\lambda)^{-\frac{\gamma}{n}}= j_N^\lambda\left((\log
\tilde h_N)'(\lambda)\right),\ \forall \lambda>0.
\end{equation}
Thus, by (\ref{ellent}), (\ref{w-elsosegben}) and
(\ref{h-elsosegben}), it turns out that
$$j_N^\lambda\left((\log
\tilde h_N)'(\lambda)\right)=\tilde
h_N(\lambda)^{-\frac{\gamma}{n}}\leq
w_N(\lambda)^{-\frac{\gamma}{n}}\leq j_N^\lambda\left((\log
w_N)'(\lambda)\right),\ \forall \lambda\in [\lambda^*,\lambda^\#].$$
Since the inverse of $j_N^\lambda$ is also increasing, it follows
that $(\log \tilde h_N)'(\lambda)\leq (\log w_N)'(\lambda)$ for
every $\lambda\in [\lambda^*,\lambda^\#]$. Therefore, the function
$\lambda\mapsto \log \frac{\tilde h_N(\lambda)}{w_N(\lambda)}$ is
non-increasing in the interval $[\lambda^*,\lambda^\#]$. In
particular, it follows that
$$0<\log \frac{\tilde h_N(\lambda^\#)}{w_N(\lambda^\#)}\leq \log \frac{\tilde h_N(\lambda^*)}{w_N(\lambda^*)}=0,$$
a contradiction, which proves the validity of the claim
(\ref{g-h-global-GN2}).

{\sc Step 6.} We shall prove that
\begin{equation}\label{asymptotic-GN2}
    \limsup_{\rho\to \infty}\frac{\textsf{m}(
    B_{x_0}(\rho))}{\omega_n\rho^n}\geq \left(\frac{\mathcal
N_{\alpha,p,n}}{\mathcal C}\right)^\frac{n}{\gamma}.
\end{equation}
By contradiction, we assume that there exists $\varepsilon_0>0$ such
that for some $\rho_0>0$,
$$\frac{\textsf{m}(
    B_{x_0}(\rho))}{\omega_n\rho^n}\leq \left(\frac{\mathcal
N_{\alpha,p,n}}{\mathcal C}\right)^\frac{n}{\gamma}-\varepsilon_0,\
\forall \rho\geq \rho_0.$$ The above inequality and
(\ref{g-h-global-GN2}) imply that for every $\lambda>\rho_0$,
\begin{eqnarray*}
  0 &\leq & w_N(\lambda)-\left(\frac{\mathcal N_{\alpha,p,n}}{\mathcal
C}\right)^\frac{n}{\gamma}h_N(\lambda) = \int_0^\lambda
\left(\frac{{\sf m}(B_{x_0}(\rho))}{\omega_n\rho^n}  -
\left(\frac{\mathcal N_{\alpha,p,n}}{\mathcal
C}\right)^\frac{n}{\gamma} \right)\omega_n\rho^n f_N(\lambda,\rho)d\rho \\
   &\leq&  \left(1+\varepsilon_0- \left(\frac{\mathcal N_{\alpha,p,n}}{\mathcal
C}\right)^\frac{n}{\gamma} \right)\int_0^{\rho_0}
\omega_n\rho^nf_N(\lambda,\rho)d\rho-\varepsilon_0 \int_0^\lambda
\omega_n\rho^nf_N(\lambda,\rho)d\rho.
\end{eqnarray*}
Reorganizing the latter estimate, it follows that for every
$\lambda>0,$
$$\varepsilon_0\frac{n}{p'}{\sf
B}\left(\frac{\alpha(p-1)+1}{1-\alpha}+1,\frac{n}{p'}\right)
\lambda^{n+p'}\leq\frac{p'}{n+p'}\left(1+\varepsilon_0-
\left(\frac{\mathcal N_{\alpha,p,n}}{\mathcal
C}\right)^\frac{n}{\gamma}
\right)\frac{\alpha(p-1)+1}{1-\alpha}\rho_0^{n+p'}.$$ Once we let
$\lambda\to \infty$, we get a contradiction. Therefore,
(\ref{asymptotic-GN2}) holds and Lemma
\ref{lemma-independence-of-point} yields that
$$\frac{\textsf{m}(
    B_{x}(\rho))}{\omega_n\rho^n}\geq \left(\frac{\mathcal
N_{\alpha,p,n}}{\mathcal C}\right)^\frac{n}{\gamma},\ \forall x\in
M,\ \rho>0,$$ which concludes the proof of Theorem {\rm
\ref{main-theorem-Gagliardo} (ii).  \hfill $\square$\\

\subsection{Limit case I $(\alpha\to 1)$:  $L^p-$logarithmic Sobolev inequality
}
In this subsection we shall provide the proof of Theorem
\ref{main-theorem-log-Sobolev}. We shall assume that $\mathcal
C>\mathcal L_{p,n}$ in $({\bf LS})_{\mathcal C}^p$.

{\sc Step 1.} As in the previous proofs, we obtain that $K=0$; the
only difference is that we shall consider
$u(x)={\textsf{m}(M)^{-{1}/{p}}}$ as a test function in $({\bf
LS})_{\mathcal C}^p$, in order to fulfil the normalization
assumption $\|u\|_{L^p}=1.$

{\sc Step 2.} Since the functions $l_p^\lambda$ ($\lambda>0$) in
Theorem B are extremals in (\ref{LS-p-generic-Dolbeau}), once we
plug them we obtain a first order ODE of the form
\begin{equation}\label{euklidesz-azonossag}
    -\log
    h_L(\lambda)+\lambda\frac{h_L'(\lambda)}{h_L(\lambda)}=\frac{n}{p}\log\left(-\mathcal
    L_{p,n}\left(\frac{p'}{p}\right)^p\lambda^p\frac{h_L'(\lambda)}{h_L(\lambda)}\right),\
    \lambda>0,
\end{equation}
where $h_L:(0,\infty)\to \mathbb R$ is defined by
$$
h_L(\lambda)=\int_{\mathbb R^n}e^{-\lambda|x|^{p'}}dx.
$$
For later use, we recall that $h_L$ can be represented alternatively
by
\begin{equation}\label{g-maskepp}
   h_L(\lambda)
   =\frac{2\pi^\frac{n}{2}}{p'\lambda^\frac{n}{p'}}\cdot\frac{\Gamma\left(\frac{n}{p'}\right)}{\Gamma\left(\frac{n}{2}\right)}=\lambda
   p'\omega_n\int_0^\infty
   e^{-\lambda\rho^{p'}}\rho^{n+p'-1}d\rho=\lambda^{-\frac{n}{p'}}p'\omega_n\int_0^\infty
  e^{-t^{p'}}t^{n+p'-1}dt.
\end{equation}

 {\sc Step 3.}  Let $w_L:(0,\infty)\to \mathbb R$ be defined by
$$
w_L(\lambda)=\displaystyle\int_M e^{-{\lambda}
\textsf{{d}}(x_0,x)^{p'}}d\textsf{m}(x),
$$
where $x_0\in M$ is the element from hypothesis {\rm $({\bf
D})^n_{x_0}$}. Note that $w_L$ is well-defined, positive and
differentiable. Indeed, by the layer cake representation, for every
$\lambda>0$ we obtain that
\begin{eqnarray*}
  w_L(\lambda) &=& \int_0^\infty \textsf{m}\left(\left\{x\in M:e^{-{\lambda}
\textsf{{d}}(x_0,x)^{p'}}>t\right\}\right)dt=\int_0^1
\textsf{m}\left(\left\{x\in M:e^{-{\lambda}
\textsf{{d}}(x_0,x)^{p'}}>t\right\}\right)dt \\
   &=&\lambda
   p'\int_0^\infty\textsf{m}(B_{x_0}(\rho))e^{-\lambda\rho^{p'}}\rho^{p'-1}d\rho\
   \ \ \ \ \ \ \ \ \ \ \ \ \ \ \ \ \ \ \ \ \ \ \  \ \ \ \ \ \ \ \ \ { [\rm change} \
   t=e^{-\lambda\rho^{p'}}]\\&\leq&
   \lambda
   p'\omega_n\int_0^\infty
   e^{-\lambda\rho^{p'}}\rho^{n+p'-1}d\rho\ \ \ \ \ \ \ \ \ \ \ \ \ \ \ \ \ \ \ \ \ \ \  \ \ \ \ \ \ \ \ \ \ \  \ \ \ \ \ \ { [\rm see\ (\ref{x_0koruli_novekedes})]}\\
   &=& h_L(\lambda)<+\infty.
\end{eqnarray*}

Let us consider the family of functions $\tilde u_\lambda:M\to
\mathbb R$ $(\lambda>0)$ defined by
$$\tilde
u_\lambda(x)=\frac{e^{-\frac{\lambda}{p}
\textsf{{d}}(x_0,x)^{p'}}}{w_L(\lambda)^\frac{1}{p}},\ x\in M.
$$
It is clear that $\| \tilde u_\lambda\|_{L^p}=1$ and as in the proof
of Theorem \ref{main-theorem-Gagliardo} (i), the function $\tilde
u_\lambda$ can be approximated by elements from Lip$_0(M);$ in fact,
$\tilde u_\lambda$ can be used as a test function in $({\bf
LS})_{\mathcal C}^p$. Thus, plugging $\tilde u_\lambda$ into the
inequality $({\bf LS})_{\mathcal C}^p,$ applying the non-smooth
chain rule and the fact that $|\nabla
\textsf{d}(x_0,\cdot)|_\textsf{d}(x)\leq 1$ for every $x\in M$, it
yields
\begin{equation}\label{metrikus-egyenlotlenseg}
    -\log
    w_L(\lambda)+\lambda\frac{w_L'(\lambda)}{w_L(\lambda)}\leq\frac{n}{p}\log\left(-\mathcal
    C\left(\frac{p'}{p}\right)^p\lambda^p\frac{w_L'(\lambda)}{w_L(\lambda)}\right),\
    \lambda>0.
\end{equation}

{\sc Step 4.} We prove that
\begin{equation}\label{egy99}
   \lim_{\lambda\to
+\infty}\frac{w_L(\lambda)}{h_L(\lambda)}= 1.
\end{equation}
For a fixed $\varepsilon>0$, let $\rho_\varepsilon>0$ from
(\ref{local-becsles-gomb}).  Then one has
\begin{eqnarray*}
  w_L(\lambda) &=& \lambda
   p'\int_0^\infty\textsf{m}(B_{x_0}(\rho))e^{-\lambda\rho^{p'}}\rho^{p'-1}d\rho
  \geq  \lambda
   p'(1-\varepsilon)\omega_n\int_0^{\rho_\varepsilon} e^{-\lambda\rho^{p'}}\rho^{n+p'-1}d\rho \\
  &=& \lambda^{-\frac{n}{p'}}p'(1-\varepsilon)\omega_n\int_0^{\rho_\varepsilon \lambda^{\frac{1}{p'}}}
  e^{-t^{p'}}t^{n+p'-1}dt.\ \ \ \ \ \ \ \ \ \ \ \ \ \ \ \ \ \ \ \ \ \ \   { [\rm change} \
   t=\lambda^\frac{1}{p'}\rho]
\end{eqnarray*}
Therefore, by the third representation of  $h_L$ (see
(\ref{g-maskepp})) it turns out that
$$\liminf_{\lambda\to
+\infty}\frac{w_L(\lambda)}{h_L(\lambda)}\geq 1-\varepsilon.$$ The
arbitrariness of  $\varepsilon>0$ together with Step 3 implies the
validity of (\ref{egy99}).

{\sc Step 5.}  We claim that
\begin{equation}\label{g-h-global-L}
    w_L(\lambda)\geq \left(\frac{\mathcal L_{p,n}}{\mathcal
C}\right)^\frac{n}{p}h_L(\lambda)=:\tilde h_L(\lambda),\ \lambda>0.
\end{equation}
Since $\mathcal C>\mathcal L_{p,n},$ by (\ref{egy99}) it follows
that  $$\lim_{\lambda\to +\infty}\frac{w_L(\lambda)}{\tilde
h_L(\lambda)}= \left(\frac{\mathcal C}{\mathcal
L_{p,n}}\right)^\frac{n}{p}>1.$$ Consequently, there exists $\tilde
\lambda>0$ such that $w_L(\lambda)>\tilde h_L(\lambda)$ for all
$\lambda>\tilde \lambda.$ If we introduce the notations
$$W(\lambda)=\log w_L(\lambda)\ {\rm and}\ \tilde H(\lambda)=\log \tilde h_L(\lambda),\ \lambda>0,$$
the latter relation implies that
\begin{equation}\label{H>G}
W(\lambda)>\tilde H(\lambda),\ \forall \lambda>\tilde \lambda,
\end{equation}
while relations in (\ref{metrikus-egyenlotlenseg}) and
(\ref{euklidesz-azonossag})
 can be rewritten in terms of $W$ and
$\tilde H$ as
\begin{equation}\label{metrikus-egyenlotlenseg-masik}
    -W(\lambda)+\lambda W'(\lambda)\leq\frac{n}{p}\log\left(-\mathcal
    C\left(\frac{p'}{p}\right)^p\lambda^p W'(\lambda)\right),\
    \lambda>0,
\end{equation}
and
\begin{equation}\label{euklidesz-azonossag-masik}
    -\tilde H(\lambda)+\lambda \tilde H'(\lambda)=\frac{n}{p}\log\left(-\mathcal
    C\left(\frac{p'}{p}\right)^p\lambda^p\tilde H'(\lambda)\right),\
    \lambda>0.
\end{equation}

Claim (\ref{g-h-global-L}) is proved once we show that
$W(\lambda)\geq \tilde H(\lambda)$ for all $\lambda>0$. By
contradiction, we assume there exists $\lambda^\#>0$ such that
$W(\lambda^\#)< \tilde H(\lambda^\#)$. Due to (\ref{H>G}),
$\lambda^\#<\tilde \lambda.$ On the one hand, let
$\lambda^*=\inf\{\lambda>\lambda^\#:W(\lambda)=\tilde H(\lambda)\}$.
 In particular,
\begin{equation}\label{ezismeg}
    W(\lambda)\leq \tilde H(\lambda),\ \forall \lambda\in
    [\lambda^\#,\lambda^*].
\end{equation}
On the other hand, if we introduce for every $\lambda>0$ the
function $j_L^\lambda:(0,\infty)\to \mathbb R$ by
$$j_L^\lambda(t)=\frac{n}{p}\log\left(\mathcal
    C\left(\frac{p'}{p}\right)^p\lambda^p t\right)+\lambda t,\ t>0,$$
relations (\ref{metrikus-egyenlotlenseg-masik}) and
(\ref{euklidesz-azonossag-masik}) become
$$-W(\lambda)\leq j_L^\lambda(-W'(\lambda))\ {\rm and}\ -\tilde H(\lambda)= j_L^\lambda(-\tilde H'(\lambda)),\ \lambda>0.$$
 By the above relations and (\ref{ezismeg}) it yields that
    $$j_L^{\lambda}(-\tilde H'(\lambda))=-\tilde H(\lambda)\leq -W(\lambda) \leq j_L^{\lambda}(-W'(\lambda)),\ \forall \lambda\in
    [\lambda^\#,\lambda^*].$$
Since $j_L^\lambda$ is increasing, it follows that $W-\tilde H$ is a
non-increasing function on $ [\lambda^\#,\lambda^*],$ which implies
$$0=(W-\tilde H)(\lambda^*)\leq (W-\tilde H)(\lambda^\#)<0,$$ a
contradiction. This completes the proof of (\ref{g-h-global-L}).

 {\sc Step 6.}  We claim
that
\begin{equation}\label{asymptotic-L}
   \limsup_{\rho\to \infty}\frac{\textsf{m}(
    B_{x_0}(\rho))}{\omega_n\rho^n}\geq \left(\frac{\mathcal L_{p,n}}{\mathcal
C}\right)^\frac{n}{p}.
\end{equation}
By assuming the contrary, there exists $\varepsilon_0>0$ such that
for some $\rho_0>0$,
$$\frac{\textsf{m}(
    B_{x_0}(\rho))}{\omega_n\rho^n}\leq \left(\frac{\mathcal L_{p,n}}{\mathcal
C}\right)^\frac{n}{p}-\varepsilon_0,\ \forall \rho\geq \rho_0.$$
Combining the latter relation with (\ref{g-h-global-L}) and
(\ref{g-maskepp}), for every $\lambda>0$ we obtain that
\begin{eqnarray*}
  0 &\leq&  w_L(\lambda)- \left(\frac{\mathcal L_{p,n}}{\mathcal
C}\right)^\frac{n}{p}h_L(\lambda) \\
   &\leq& \lambda
   p' \int_0^{\rho_0}
   \textsf{m}(B_{x_0}(\rho))e^{-\lambda\rho^{p'}}\rho^{p'-1}d\rho
   +\lambda
   p'\omega_n \left(\left(\frac{\mathcal L_{p,n}}{\mathcal
C}\right)^\frac{n}{p}-\varepsilon_0\right)\int_{\rho_0}^\infty
   e^{-\lambda\rho^{p'}}\rho^{n+p'-1}d\rho\\&&-
   \lambda
   p'\omega_n\left(\frac{\mathcal L_{p,n}}{\mathcal
C}\right)^\frac{n}{p}\int_0^\infty
   e^{-\lambda\rho^{p'}}\rho^{n+p'-1}d\rho.
\end{eqnarray*}
Rearranging the above inequality, by virtue of
(\ref{x_0koruli_novekedes}) it follows for every $\lambda>0$ that
$$\varepsilon_0\int_0^\infty
   e^{-\lambda\rho^{p'}}\rho^{n+p'-1}d\rho\leq\left(1-\left(\frac{\mathcal L_{p,n}}{\mathcal
C}\right)^\frac{n}{p}+\varepsilon_0\right)\int_0^{\rho_0}
   e^{-\lambda\rho^{p'}}\rho^{n+p'-1}d\rho.$$
   Due to  (\ref{g-maskepp}), the latter inequality implies
   $$\varepsilon_0\frac{1}{p'\lambda^{1+\frac{n}{p'}}}\Gamma\left(\frac{n}{p'}+1\right)\leq\left(1-\left(\frac{\mathcal L_{p,n}}{\mathcal
C}\right)^\frac{n}{p}+\varepsilon_0\right)\frac{\rho_0^{n+p'}}{n+p'},\
\lambda>0.$$ Now, letting $\lambda\to 0^+$ we arrive to a
contradiction. Therefore, the proof of (\ref{asymptotic-L})  is
concluded. Thus, Lemma \ref{lemma-independence-of-point} gives that
$$\frac{\textsf{m}(
    B_{x}(\rho))}{\omega_n\rho^n}\geq \left(\frac{\mathcal
L_{p,n}}{\mathcal C}\right)^\frac{n}{p},\ \forall x\in M,\ \rho>0,$$
concluding the proof of Theorem \ref{main-theorem-log-Sobolev}.
 \hfill $\square$\\

\subsection{Limit case II $(\alpha\to 0)$: Faber-Krahn-type inequality
} In this part we sketch the proof of
Theorem \ref{main-theorem-Faber-Krahn}. Similarly as before, we
assume that $\mathcal C>\mathcal F_{p,n}.$

 {\sc Step 1.} Analogously to
Theorem \ref{main-theorem-Gagliardo} (i), it follows that $K=0.$

{\sc Step 2.}  The function $x\mapsto
\left(\lambda^{p'}-|x|^{p'}\right)_+$ being extremal in
(\ref{Faber-Krahn}) for every $\lambda>0$, a direct computation
shows that
\begin{eqnarray}\label{1111-FK}
 h_F(\lambda)  &=& \mathcal F_{p,n}p'\left(-h_F(\lambda)+\frac{1}{p'}\lambda h_F'(\lambda)\right)^\frac{1}{p} \left(\frac{1}{p'}
\lambda^{1-p'}h_F'(\lambda)\right)^{1-\frac{1}{p^\star}},
\end{eqnarray}
where  $h_F:(0,\infty)\to \mathbb R$ is given by
$$h_F(\lambda)=\int_{\mathbb R^n}\left(\lambda^{p'}-|x|^{p'}\right)_+dx,\
\lambda>0.$$

{\sc Step 3.} Let  $x_0\in M$ from {\rm $({\bf D})^n_{x_0}$}.  Since
$u_\lambda= \left(\lambda^{p'}-{\sf d}(x_0,\cdot)^{p'}\right)_+\in
{\rm Lip}_0(M)$, we may insert $u_\lambda$ into ${\bf
(FK)}_{\mathcal C}^{p}$ obtaining
\begin{equation}\label{fk-egyenlotlenseg}
    \|u_\lambda\|_{L^1}\leq \mathcal C\||\nabla
u_\lambda|_{\sf{{d}}}\|_{L^p}{\sf m}({\rm
supp}(u_\lambda))^{1-\frac{1}{p^\star}}.
\end{equation}
 First, we observe that
 $$|\nabla u_\lambda|_{\sf{{d}}}(x)= p'{\sf
d}(x_0,x)^{p'-1}|\nabla {\sf d}(x_0,\cdot)|_{\sf{{d}}}(x)\leq p'{\sf
d}(x_0,x)^{p'-1},\ \forall x\in B_{x_0}(\lambda),$$ while $|\nabla
u_\lambda|_{\sf{{d}}}(x)=0$ for every $x\notin B_{x_0}(\lambda)$.
Moreover, since the spheres have zero ${\sf m}-$measures (see
Theorem \ref{CD-tetel}), we have that
$${\sf m}({\rm
supp}(u_\lambda))={\sf m}(\overline{B_{x_0}(\lambda)})={\sf
m}({B_{x_0}(\lambda)}).$$ We now introduce the function
$w_F:(0,\infty)\to \mathbb R$  given by
$$w_F(\lambda)=\int_{M}\left(\lambda^{p'}-{\sf d}(x_0,x)^{p'}\right)_+d{\sf m}(x),\
\lambda>0.$$ Due to the layer cake representation, one has
\begin{eqnarray*}
 w_F(\lambda)  &=&  \int_{B_{x_0}(\lambda)}\left(\lambda^{p'}-{\sf d}(x_0,x)^{p'}\right)d{\sf m}(x)=\lambda^{p'}{\sf
m}({B_{x_0}(\lambda)})-\int_{B_{x_0}(\lambda)}{\sf d}(x_0,x)^{p'}d{\sf m}(x)\\
   &=& \lambda^{p'}{\sf
m}({B_{x_0}(\lambda)})- \int_0^{\lambda^{p'}} {\sf m}\left(\{x\in B_{x_0}(\lambda):{\sf d}(x_0,x)^{p'}>t\}\right)dt\\
   &=& \lambda^{p'}{\sf
m}({B_{x_0}(\lambda)})- p'\int_0^{\lambda} \left({\sf
m}({B_{x_0}(\lambda)})-{\sf m}({B_{x_0}(\rho)})\right) \rho^{p'-1}
d\rho\ \ \ \ \ [{\rm change}\ t=\rho^{p'}] \\&=&p'\int_0^\lambda
{\sf m}({B_{x_0}(\rho)}) \rho^{p'-1} d\rho.
\end{eqnarray*}
Therefore,
$$\|u_\lambda\|_{L^1}=w_F(\lambda),\ {\sf m}({\rm
supp}(u_\lambda))={\sf
m}({B_{x_0}(\lambda)})=\frac{1}{p'}\lambda^{1-p'}w_F'(\lambda),$$
and
$$\||\nabla u_\lambda|_{\sf{{d}}}\|_{L^p}\leq p'\left(\int_{B_{x_0}(\lambda)}{\sf d}(x_0,x)^{p'}d{\sf m}(x)\right)^\frac{1}{p}=p'\left(-w_F(\lambda)+\frac{1}{p'}\lambda w_F'(\lambda)\right)^\frac{1}{p}.$$
 Consequently,  inequality
(\ref{fk-egyenlotlenseg}) takes the form
$$w_F(\lambda)  \leq \mathcal Cp'\left(-w_F(\lambda)+\frac{1}{p'}\lambda w_F'(\lambda)\right)^\frac{1}{p} \left(\frac{1}{p'}
\lambda^{1-p'}w_F'(\lambda)\right)^{1-\frac{1}{p^\star}},\
\lambda>0,$$ which is formally (\ref{2222-GN2}) if $\alpha\to 0$
since due to (\ref{gamma-best}), $\lim_{\alpha\to 0}\gamma=1$ and
$\lim_{\alpha\to 0}\frac{1-\gamma}{\alpha p}=1-\frac{1}{p^\star}.$

Therefore,  we may proceed as in the proof of Theorem
\ref{main-theorem-Gagliardo} (ii) (Steps 4-6), proving that
$$ \lim_{\lambda\to 0^+}\frac{w_F(\lambda)}{h_F(\lambda)}=1,$$
$$ w_F(\lambda)\geq \left(\frac{\mathcal F_{p,n}}{\mathcal
C}\right)^n h_F(\lambda),\ \forall \lambda>0,$$ and finally
$$\frac{\textsf{m}(
    B_{x}(\rho))}{\omega_n\rho^n}\geq \left(\frac{\mathcal
F_{p,n}}{\mathcal C}\right)^n,\ \forall x\in M,\ \rho>0,$$ which
concludes the proof of Theorem \ref{main-theorem-Faber-Krahn}.
\hfill $\square$

\section{\sc Rigidity results in smooth settings}\label{sect-4}

As a starting point, we need  an Aubin-Hebey-type result (see
\cite{Aubin} and \cite{Hebey}) for Gagliardo-Nirenberg inequalities
which is valid on generic Riemannian manifolds.

\begin{lemma}\label{lemma-Riemannian}
Let $(M,g)$ be a complete $n-$dimensional Riemannian manifold and
$\mathcal C>0$. The following statements hold:
\begin{itemize}
\item[{\rm (i)}] If ${\bf (GN1)}_{\mathcal C}^{\alpha, p}$ holds on $(M,g)$ for some $p\in (1,n)$ and $\alpha\in
(1,\frac{n}{n-p}]$ then $\mathcal C\geq \mathcal G_{\alpha,p,n};$
\item[{\rm (ii)}] If  ${\bf (GN2)}_{\mathcal C}^{\alpha, p}$ holds on $(M,g)$ for some $p\in (1,n)$ and $\alpha\in
(0,1)$  then $\mathcal C\geq \mathcal N_{\alpha,p,n};$
  \item[{\rm (iii)}] If $({\bf
LS})_{\mathcal C}^p$ holds on $(M,g)$ for some $p\in (1,n)$  then
$\mathcal C\geq \mathcal L_{p,n};$
 \item[{\rm (iv)}] If $({\bf
FK})_{\mathcal C}^p$ holds on $(M,g)$ for some $p\in (1,n)$  then
$\mathcal C\geq \mathcal F_{p,n}.$
\end{itemize}
\end{lemma}

{\it Proof.} (i) By contradiction, we assume that ${\bf
(GN1)}_{\mathcal C}^{\alpha, p}$ holds on $(M,g)$ for some $p\in
(1,n)$, $\alpha\in (1,\frac{n}{n-p}]$,  and $\mathcal C< \mathcal
G_{\alpha,p,n}.$ Let $x_0\in M$ be fixed arbitrarily. For every
$\varepsilon>0$, there exists a local chart $(\Omega,\phi)$ of $M$
at the point $x_0$ and a number $\delta>0$ such that
$\phi(\Omega)=B_0(\delta)=\{\tilde x\in \mathbb R^n:|\tilde
x|<\delta\}$ and the components $g_{ij}=g_{ij}(x)$ of the Riemannian
metric $g$ on  $(\Omega,\phi)$ satisfy
\begin{equation}\label{two-sided}
    (1-\varepsilon)\delta_{ij}\leq g_{ij} \leq (1+\varepsilon)\delta_{ij}
\end{equation}
in the sense of bilinear forms. Since ${\bf (GN1)}_{\mathcal
C}^{\alpha, p}$ is valid, relation (\ref{two-sided}) shows that for
every $\varepsilon>0$ small enough, there exists
$\delta_\varepsilon>0$ and $\mathcal C_\varepsilon\in (\mathcal C,
\mathcal G_{\alpha,p,n})$ such that for every $\delta\in
(0,\delta_\varepsilon)$ and $v\in {\rm Lip}_0(B_0(\delta))$,
\begin{equation}\label{c-hpw-uj-meg}
      \|v\|_{L^{\alpha p}(B_0(\delta),dx)}\leq \mathcal C_\varepsilon
    \|\nabla
    v\|_{L^p(B_0(\delta),dx)}^{\theta}\|v\|_{L^{\alpha(p-1)+1}(B_0(\delta),dx)
    }^{1-\theta}.
\end{equation}
Let us fix $u\in {\rm Lip}_0(\mathbb R^n)$  arbitrarily and set
$v_\lambda(x)=\lambda^\frac{n}{p} u(\lambda x)$, $\lambda>0.$ For
$\lambda>0$ large enough, one has $v_\lambda\in {\rm
Lip}_0(B_0(\delta))$. If we plug in $v_\lambda$ into
(\ref{c-hpw-uj-meg}), by using the scaling properties
\begin{equation}\label{scaling-elol}
    \|\nabla
v_\lambda\|_{L^p(B_0(\delta),dx)}=\lambda\|\nabla u\|_{L^p(\mathbb
R^n,dx)}\ \ {\rm and}\ \  \|
v_\lambda\|_{L^q(B_0(\delta),dx)}=\lambda^{\frac{n}{p}-\frac{n}{q}}
\| u\|_{L^q(\mathbb R^n,dx)},\ \forall q>0,
\end{equation}
 and the form of the number $\theta$ (see (\ref{theta-best})), it
follows that
$$ \|u\|_{L^{\alpha p}(\mathbb
R^n,dx)}\leq \mathcal C_\varepsilon
    \|\nabla
    u\|_{L^p(\mathbb
R^n,dx)}^{\theta}\|u\|_{L^{\alpha(p-1)+1}
    (\mathbb
R^n,dx)}^{1-\theta}.$$
 If we insert the extremal function $h_{\alpha,p}^\lambda$ of the optimal Gagliardo-Nirenberg inequality on $\mathbb R^n$ ($\alpha>1$) into the latter relation,  Theorem A  yields
 that $\mathcal G_{\alpha,p,n}\leq \mathcal C_\varepsilon$, a
 contradiction.

 The proofs of (ii) (iii) and (iv) are analogous to (i), taking into
 account  in addition to (\ref{scaling-elol})
  that
 $${\bf Ent}_{dx}(|v_\lambda|^p)={\bf Ent}_{dx}(|u|^p)+n\|u\|_{L^p}^p\log \lambda,$$
 and
 $$\mathcal H^n({\rm supp}(v_\lambda))=\lambda^{-n}\mathcal
H^n({\rm supp}(u)),$$ respectively.
  \hfill $\square$\\

\subsection{Gagliardo-Nirenberg inequalities on Riemannian manifolds with Ricci$\geq
0$}\label{sect-riemann-3} Before presenting the proofs of Theorem
\ref{theorem-Riemann-log-Sobolev} and Corollary
\ref{corollary-Riemann}, we recall some results from Munn
\cite{Munn-JGA}.

To do this, let $(M,g)$ be an $n(\geq 2)-$dimensional complete
Riemannian manifold with non-positive Ricci curvature endowed with
its canonical volume element $dv_g$. The {\it asymptotic volume
growth} of $(M,g)$ is defined by
$${\rm AVG}_{(M,g)}=\lim_{r\to \infty}\frac{{\rm Vol}_g(B_x(r))}{\omega_n r^n}.$$
By Bishop-Gromov comparison theorem it follows that ${\rm
AVG}_{(M,g)}\leq 1$ and this number is independent of the point
$x\in M.$

Given $k\in\{1,...,n\}$, let us denote by $\delta_{k,n}>0$  the
smallest positive solution to the equation
$10^{k+2}C_{k,n}(k)s\left(1+\frac{s}{2k}\right)^k=1$ in variable
$s$,  where
\begin{equation*}
C_{k,n}(i)= \left\{
\begin{array}{lll}
1 & \mbox{if} & i=0, \\
3+10C_{k,n}(i-1)+(16k)^{n-1}(1+10C_{k,n}(i-1))^n & \mbox{if} & i\in
\{1,...,k\}.
\end{array}%
\right.
\end{equation*} We now consider
the smooth, bijective and increasing function
$h_{k,n}:(0,\delta_{k,n})\to (1,\infty)$ defined by
$$h_{k,n}(s)=\left[1-10^{k+2}C_{k,n}(k)s\left(1+\frac{s}{2k}\right)^k\right]^{-1}.$$
For every $s>1,$ let
\begin{equation*}
\beta(k, s,n)= \left\{
\begin{array}{lll}
1-\left[1+\frac{s^n}{[h_{1,n}^{-1}(s)]^n}\right]^{-1} & \mbox{if} & k=1, \\
\max\left\{\beta(1, s,n),\beta(i,
1+\frac{h_{k,n}^{-1}(s)}{2k},n):i=1,...,k-1\right\} & \mbox{if} &
k\in \{2,...,n\}.
\end{array}%
\right.
\end{equation*}
Note that the constant $\beta(k, s,n)$,  which is used to prove the
Perelman's maximal volume lemma, denotes the minimum volume growth
of $(M,g)$ needed to guarantee that any continuous map $f : \mathbb
S^k \to B_x(\rho)$ has a continuous extension $g : \mathbb D
^{k+1}\to B_x(c\rho),$ where $\mathbb D ^{k+1}=\{y\in \mathbb
R^{k+1}:|y|\leq 1\}$ and $\mathbb S^k=\partial \mathbb D ^{k+1},$
see \cite[Definition 3.3]{Munn-JGA}.  Finally, the {\it
Munn-Perelman constant} is defined as
$$\alpha_{MP}(k,n)=\inf_{s\in (1,\infty)}\beta(k, s,n).$$
By construction, $\alpha_{MP}(k,n)$ is non-decreasing in $k$; for
numerical values of $\alpha_{MP}(k,n)$ one can consult
\cite[Appendix A]{Munn-JGA}.\\

{\it Proof of Theorem \ref{theorem-Riemann-log-Sobolev}.} Let
$(M,g)$ be an $n-$dimensional complete Riemannian manifold
 with non-negative Ricci
curvature $(n\geq 2)$ and assume the $L^p-$logarithmic Sobolev
inequality $({\bf LS})_{\mathcal C}^p$ holds on $(M,g)$ for some
$p\in (1,n)$ and $\mathcal C>0$.

(i) It follows from Lemma \ref{lemma-Riemannian} (iii), i.e.,
$\mathcal C\geq \mathcal L_{p,n}.$

(ii) Anderson \cite{Anderson} and Li \cite{Peter_Li} stated that if
there exists $c_0>0$ such that Vol$_g(B_x(\rho))\geq c_0 \omega_n
\rho^n$ for every $\rho>0$, then $(M,g)$ has finite fundamental
group $\pi_1(M)$ and its order is bounded above by ${c_0}^{-1}.$
Thus it remains to apply Theorem \ref{main-theorem-log-Sobolev}.

(iii) Assume that $\mathcal C <
\alpha_{MP}(k_0,n)^{-\frac{p}{n}}\mathcal L_{p,n}$ for some $k_0\in
\{1,...,n\}$. By Theorem \ref{main-theorem-log-Sobolev}, we have
that
$${\rm AVG}_{(M,g)}=\lim_{r\to \infty}\frac{{\rm Vol}_g(B_x(r))}{\omega_n r^n}\geq\left(\frac{\mathcal L_{p,n}}{\mathcal C}\right)^\frac{n}{p}>\alpha_{MP}(k_0,n)\geq...\geq \alpha_{MP}(1,n).$$
By Munn \cite[Theorem 1.2]{Munn-JGA}, it follows that
$\pi_1(M)=...=\pi_{k_0}(M)=0.$

(iv) If $\mathcal C < \alpha_{MP}(n,n)^{-\frac{p}{n}}\mathcal
L_{p,n}$, then $\pi_1(M)=...=\pi_{n}(M)=0,$ which implies the
contractibility of $M$, see e.g. Luft \cite{Luft}.

(v) If $\mathcal C = \mathcal L_{p,n}$ then by Theorem
\ref{main-theorem-log-Sobolev} and the Bishop-Gromov volume
comparison theorem follows that Vol$_g(B_x(\rho))=\omega_n\rho^n$
for every $x\in M$ and $\rho>0$. The equality in Bishop-Gromov
theorem implies that $(M,g)$ is isometric to the Euclidean space
$\mathbb R^n$. The converse
trivially holds. \hfill $\square$\\

\begin{remark}\rm\label{rem-3.1}
 In the study of heat kernel bounds on an $n-$dimensional complete Riemannian manifold $(M,g)$
 with non-negative Ricci
curvature,  the logarithmic Sobolev inequality
\begin{equation}\label{LS-p=2}
    {\bf Ent}_{dv_g}(u^2)\leq \frac{n}{2}\log\left(C\|\nabla_g u\|^2_{L^2(M,dv_g)}\right),\ \forall u\in C_0^\infty(M),\
\|u\|_{L^2}=1,
\end{equation}
plays a central role, $C>0$. In fact, (\ref{LS-p=2}) is equivalent
to an upper bound of the heat kernel $p_t(x,y)$ on $M$, i.e.,
\begin{equation}\label{heat-kernel}
    \sup_{x,y\in M}p_t(x,y)\leq \tilde C t^{-\frac{n}{2}},\ t>0,
\end{equation}
for some $\tilde C>0$. According to Theorem B (from
\S\ref{sect1.1}), the optimal constant in (\ref{LS-p=2}) in the
Euclidean space $\mathbb R^n$  is given by $C=\mathcal
L_{n,2}=\frac{2}{n\pi e };$ this scale invariant form on $\mathbb
R^n$ can be deduced by Gross \cite{Gross} logarithmic Sobolev
inequality
$${\bf Ent}_{d\gamma_n}(u^2)\leq 2\|\nabla u\|_{L^2(\mathbb R^n,d\gamma_n)}^2,\ \forall u\in C_0^\infty(\mathbb R^n),\
\|u\|_{L^2(\mathbb R^n,d\gamma_n)}=1,$$ where the canonical Gaussian
measure $\gamma_n$ has the density
$\delta_n(x)=(2\pi)^{-\frac{n}{2}}e^{-\frac{|x|^2}{2}},$ $x\in
\mathbb R^n,$ see Weissler \cite{Weissler}. Sharp estimates on the
heat kernel shows that {\it on a complete Riemannian manifold
$(M,g)$ with non-negative Ricci curvature the $L^2-$logarithmic
Sobolev inequality (\ref{LS-p=2}) holds with the optimal Euclidean
constant $C=\mathcal L_{n,2}=\frac{2}{n\pi e }$ if and only if
$(M,g)$ is isometric to $\mathbb R^n$}, cf. Bakry, Concordet and
Ledoux \cite{BCL}, Ni \cite{Ni}, and Li \cite{Peter_Li}. In this
case, $\tilde C=(4\pi)^{-\frac{n}{2}}$ in (\ref{heat-kernel}).

In particular, Theorem \ref{theorem-Riemann-log-Sobolev} (v) gives a
positive answer to the open problem of C. Xia \cite{Xia} concerning
 the validity of the optimal
$L^p-$logarithmic So\-bo\-lev
 inequality for generic $p\in (1,n)$ in the same geometric context as above.
 Xia's formulation was deeply motivated by the lack of sharp $L^p-$estimates $(p\neq
2)$ for the heat kernel on Riemannian manifolds with non-negative
Ricci curvature.
%
%
\end{remark}

Similar results to Theorem \ref{theorem-Riemann-log-Sobolev} can be
stated for the other three Gagliardo-Nirenberg-type  inequalities;
here we formulate one for ${\bf (GN1)}_{\mathcal C}^{\alpha,p}$, the
other two inequalities are left to the reader.

\begin{theorem}\label{theorem-Riemann-Gagliardo}
Let $(M,g)$ be an $n-$dimensional complete Riemannian manifold
 with non-negative Ricci
curvature $(n\geq 2)$ and assume the ${\bf (GN1)}_{\mathcal
C}^{\alpha,p}$ holds on $(M,g)$ for some $p\in (1,n)$, $\alpha\in
(1,\frac{n}{n-p}]$ and $\mathcal C>0$. Then the following assertions
hold:
\begin{itemize}
 \item[{\rm (i)}] $\mathcal C\geq \mathcal G_{\alpha,p,n};$
  \item[{\rm (ii)}] The order of the fundamental group $\pi_1(M)$ is bounded
  above by $\left(
\frac{\mathcal C}{\mathcal G_{\alpha,p,n}}
\right)^\frac{n}{\theta};$
 \item[{\rm (iii)}] If $\mathcal C < \alpha_{MP}(k_0,n)^{-\frac{\theta}{n}}\mathcal G_{\alpha,p,n}$ for some $k_0\in \{1,...,n\}$ then
 $\pi_1(M)=...=\pi_{k_0}(M)=0;$
 \item[{\rm (iv)}] If $\mathcal C < \alpha_{MP}(n,n)^{-\frac{\theta}{n}}\mathcal G_{\alpha,p,n}$  then
 $M$ is contractible$;$
  \item[{\rm (v)}]
  $\mathcal C=\mathcal G_{\alpha,p,n}$   if and only if $(M,g)$ is isometric to the Euclidean space $\mathbb
  R^n.$
\end{itemize}
\end{theorem}

%

\subsection{Gagliardo-Nirenberg inequalities on Finsler manifolds with $n-$Ricci$\geq
0$}\label{sect-finsler}

Let $M$ be a connected $n-$dimensional $C^{\infty}$-manifold and
$TM=\bigcup_{x \in M}T_{x} M $ be its tangent bundle. The pair
$(M,F)$ is called a \textit{reversible Finsler manifold} if a
continuous function $F:TM\longrightarrow [0,\infty)$ satisfies the
conditions:

(a) $F\in C^{\infty}(TM\setminus \{0\})$;

(b) $F(x,tv)=|t|F(x,v)$ for all $t\in \mathbb R$ and $(x,v)\in TM$;

(c) the $n \times n$ matrix $
g_{ij}(x,v)=\frac{1}{2}\frac{\partial^2 (F^2)}{\partial v^{i}
\partial v^{j}}(x,v)$ is positive definite for all $(x,v)\in TM\setminus\{ 0 \}$.

\noindent Here $ v=\sum_{i=1}^n v^i \frac{\partial}{\partial x^i}, $
and we shall denote by $g_v$ the inner product on $T_xM$ induced by
the above form. If $g_{ij}(x)=g_{ij}(x,v)$ is independent of $v$
then $(M,F)$ is called {\it Riemannian manifold}. A {\it Minkowski
space} consists of a finite dimensional vector space $V$ and a
Minkowski norm which induces a Finsler metric on $V$ by translation,
i.e., $F(x,v)$ is independent of $x$. A Finsler manifold $(M,F)$ is
called a {\it locally Minkowski space} if every point in $M$ admits
a local coordinate system $(x^i)$ on its neighborhood such that
$F(x,v)$ depends only on $v$ and not on $x$.

We consider on the pull-back bundle $\pi ^{*}TM$ the \textit{Chern
connection}, see Bao, Chern and Shen \cite[Theorem 2.4.1]{BCS}. The
coefficients of the Chern connection are denoted by
$\Gamma_{jk}^{i}$, which are instead of the well-known Christoffel
symbols from Riemannian geometry. A Finsler manifold is of
\textit{Berwald type} if the coefficients $\Gamma_{ij}^{k}(x,v)$ in
natural coordinates are independent of $v$. It is clear that
Riemannian manifolds and $($locally$)$ Minkowski spaces are Berwald
spaces. The Chern connection induces in a natural manner on $\pi
^{*}TM$ the \textit{curvature tensor} $R$, see Bao, Chern and Shen
\cite[Chapter 3]{BCS}. By means of the connection, we also have
 the {\it covariant derivative} $D_vu$ of a vector field
$u$ in the direction $v\in T_xM.$  Note that $v\mapsto D_vu$ is not
linear. A vector field $u=u(t)$ along a curve $\sigma$ is {\it
parallel} if $D_{\dot \sigma}u=0.$ A $C^\infty$ curve
$\sigma:[0,a]\to M$ is a {\it geodesic} if  $D_{\dot \sigma}{\dot
\sigma}=0.$ Geodesics are considered to be parametrized
proportionally to arc-length. The Finsler manifold is {\it complete}
if every geodesic segment can be extended to $\mathbb R.$ For a
$C^{\infty}$-curve $\sigma: [0,l]\longrightarrow M$, its {integral
length} is given by $L_F(\sigma):=\displaystyle\int_{0}^{l}
F(\sigma(t), \dot\sigma(t))\,dt$. Define the {\it distance function}
$d_{F}: M\times M \longrightarrow[0,\infty)$ by $$d_{F}(x_1,x_2) =
\inf_{\sigma} L_F(\sigma),$$ where $\sigma$ runs over all
$C^{\infty}$-curves from $x_1$ to $x_2$. Geodesics locally minimize
$d_F-$dis\-tan\-ces.

Let $u,v\in T_xM$ be two non-collinear vectors and $\mathcal S={\rm
span}\{u,v\}\subset T_xM$. By means of the curvature tensor $R$, the
{\it flag curvature} of the flag $\{\mathcal S,v\}$ is defined by
$$
K(\mathcal S;v) =\frac{g_v(R(U,V)V, U)}{g_v(V,V) g_v(U,U) -
g_v(U,V)^{2}},
$$
where $U=(v;u),V=(v;v)\in \pi^*TM.$  If $(M,F)$ is Riemannian, the
flag curvature reduces to the well known sectional curvature.

Let $v\in T_xM$ be such that $F(x,v)=1$ and let
$\{e_i\}_{i=1,...,n}$ with $e_n=v$ be a basis for $T_xM$ such that
$\{(v;e_i)\}_{i=1,...,n}$ is an orthonormal basis for $\pi_*TM$. Let
$\mathcal S_i={\rm span}\{e_i,v\}$, $i=1,...,n-1.$ The {\it Ricci
curvature} Ric$:TM\to \mathbb R$  is defined by ${\rm
Ric}(cv)=c^2\sum_{i=1}^{n-1}K(\mathcal S_i;v)$ for every $c> 0.$

Let $(M,F)$ be an $n-$dimensional complete Finsler manifold and let
$\textsf{m}$ be an arbitrarily positive smooth measure on $M$; such
a manifold is viewed as a regular metric measure space and we denote
it by $(M,F,\textsf{m})$. Let $v\in T_xM$ be such that $F(x,v)=1$
and let
$$\Upsilon(v)=\log\left(\frac{{\rm vol}_{g_v}({\rm
B}(0,1))}{\textsf{m}_x({\rm B}(0,1))}\right),$$ where ${\rm
vol}_{g_v}$ and $\textsf{m}_x$ denote the Lebesgue measures on
$T_xM$ induced by $g_v$ and $\textsf{m}$, respectively, while ${\rm
B}(0,1)=\{y\in T_xM:F(x,y)< 1\}$ is the unit tangent ball at $T_xM$.
The latter relation can be rewritten into the more familiar form
$\textsf{m}_x({\rm B}(0,1))=e^{-\Upsilon(v)}{\rm vol}_{g_v}({\rm
B}(0,1)).$ We introduce the notation
\begin{equation}\label{upsilon}
    \partial_v \Upsilon=\frac{d}{dt}\Upsilon (\dot\sigma(t))\big|_{t=0},
\end{equation}
where $\sigma:(-\varepsilon,\varepsilon)\to M$ is the geodesic with
$\sigma(0)=x$ and $\dot\sigma(0)=v$. We say that the space
$(M,F,\textsf{m})$ has $n-${\it Ricci curvature bounded below by}
$K\in \mathbb R$ if and only if ${\rm Ric}(v)\geq K$ and $\partial_v
\Upsilon=0$ for every $v\in T_xM$ such that $F(x,v)=1$, see Ohta
\cite[Theorem 1.2]{Ohta} and Ohta and Sturm \cite[Definition
5.1]{Ohta-Sturm}. Note that a Berwald space endowed with  the
Busemann-Hausdorff measure $\textsf{m}_{ BH}$ (and inducing the
volume form $dV_F$)
verifies the property $ \partial_v \Upsilon\equiv 0,$ see Shen
\cite[Propositions 2.6 \& 2.7]{Shen-volume}.

The {\it polar transform}  of $F$ is defined for every
$(x,\alpha)\in T^*M$ by
\begin{equation}\label{polar-trans}
    F^*(x,\alpha)
 =\sup_{v\in T_xM\setminus \{0\}}\frac{\alpha(v)}{F(x,v)}.
\end{equation}
Note that, for every $x\in M$, the function $F^*(x,\cdot)$ is a
Minkowski norm on $T_x^*M$.

If $u\in {\rm Lip}_0(M)$, then relation (\ref{local-constant}) can
be interpreted as
\begin{equation}\label{finsler-derivaltak}
    |\nabla u|_{d_F}(x)=F^*(x,Du(x))\ {\rm for\ a.e.}\ x\in M,
\end{equation}
 where $Du(x)\in T_x^*(M)$ is the
distributional derivative of $u$ at $x\in M$, see Ohta and Sturm
\cite{Ohta-Sturm}.  In particular, if $(M,F)=(M,g)$ is a Riemannian
manifold, then $|\nabla u|_{d_g}=|\nabla_g u|$, where $d_g$ is the
distance function on $(M,g)$,  $\nabla_g$ is the Riemannian gradient
on $(M,g)$, and $|\cdot|$ is the norm coming from the Riemannian
metric $g$, respectively. \\

Although a slightly more general result can be proved, we present an
application on  Berwald spaces $(M,F)$ endowed with the canonical
Busemann-Hausdorff measure ${\sf m}_{BH}$ (and its induced volume
form $dV_F$), by exploring the results of Cordero-Erausquin, Nazaret
and Villani \cite{CE-N-Villani} and Gentil \cite{Gentil} (see
Theorems A \& B).

\begin{theorem}\label{theorem-Finsler-1} {\rm [Optimality vs. flatness]} Let $(M,F)$ be an $n-$dimensional complete reversible Berwald space with non-negative Ricci
curvature. The following statements are equivalent:
\begin{itemize}
\item[{\rm (i)}]  ${\bf (GN1)}_{\mathcal G_{\alpha,p,n}}^{\alpha, p}$ holds on $(M,F)$ for some $p\in (1,n)$ and $\alpha\in
(1,\frac{n}{n-p}];$
\item[{\rm (ii)}]  ${\bf (GN2)}_{\mathcal N_{\alpha,p,n}}^{\alpha, p}$ holds on $(M,F)$ for some $p\in (1,n)$ and $\alpha\in
(0,1);$
  \item[{\rm (iii)}]  $({\bf
LS})_{\mathcal L_{p,n}}^p$ holds on $(M,F)$ for some $p\in (1,n);$
 \item[{\rm (iv)}]  $({\bf
FK})_{\mathcal F_{p,n}}^p$ holds on $(M,F)$ for some $p\in (1,n);$
  \item[{\rm (v)}] $(M,F)$ is isometric to an
$n-$di\-mensional Minkowski space.
\end{itemize}
\end{theorem}

{\it Proof.}  We divide the proof into two parts.

(i)$\vee$(ii)$\vee$(iii)$\vee$(iv)$\Rightarrow$(v). Note that the
Busemann-Hausdorff measure ${\sf m}_{BH}$ satisfies the
$n-$den\-sity assumption for every $x\in M$, i.e.,
$$\lim_{\rho\to
 0}\frac{\textsf{{m}}_{BH}( B_{x}(\rho))}{\omega_n\rho^n}=1,$$
see Shen \cite[Lemma 5.2]{Shen-volume}.  Since $(M,F)$ is a Berwald
space (thus $
\partial_v \Upsilon\equiv 0$ for every $v\in T_xM$, $x\in M$, see (\ref{upsilon})), the non-negativity of the Ricci curvature on $(M,F)$
coincides with the non-negativity of the $n-$Ricci curvature on
$(M,d_F,\textsf{{m}}_{BH})$, thus the metric measure space
$(M,d_F,\textsf{{m}}_{BH})$ satisfies the curvature-dimension
condition $\textsf{CD}(0,n)$, see Ohta \cite{Ohta}. Moreover, the
completeness of $(M,F)$ via Hopf-Rinow theorem implies that the
$(M,d_F,\textsf{{m}}_{BH})$ is proper. Applying now any of the
Theorems \ref{main-theorem-Gagliardo},
\ref{main-theorem-log-Sobolev} or \ref{main-theorem-Faber-Krahn}
(according to which of the assumptions (i), (ii), (iii) or (iv) is
satisfied), it yields that
$$\textsf{{m}}_{BH}( B_x(\rho))\geq \omega_n\rho^n\ {\rm for\ all}\ x\in
M,\ \rho\geq 0.$$  By the generalized Bishop-Gromov theorem on
Finsler manifolds and the $n-$density property we also have the
reverse inequality, thus
\begin{equation}\label{identitas}
   \textsf{{m}}_{BH}( B_x(\rho))=
\omega_n\rho^n\ {\rm for\ all}\ x\in M,\ \rho\geq 0.
\end{equation}
The latter relation immediately implies that the flag curvature on
$(M,F)$ is identically zero, see Ohta \cite[Theorem 7.3]{Ohta}, and
Krist\'aly and Ohta \cite[Theorem 3.3]{Kri-Ohta}. Due to Bao, Chern
and Shen \cite[Section~10.5]{BCS}), every Berwald space with zero
flag curvature is necessarily a locally Minkowski space. By
(\ref{identitas}) it follows that $(M,F)$ is actually isometric to a
Minkowski space.

(v)$\Rightarrow$(i)$\wedge$(ii)$\wedge$(iii)$\wedge$(iv).  Let us
fix an arbitrary norm $\|\cdot\|$ on $\mathbb R^n$, and let $\Phi
:(M,F)\rightarrow (\mathbb{R}^{n},\Vert \cdot \Vert )$ be an
isometry.  Then
$$
F(x,y)=\Vert d\Phi _{x}(y)\Vert ,\ x\in M,y\in T_{x}M,
$$
 and a simple computation based on
the definition of the polar transform (see (\ref{polar-trans}))
gives
\begin{equation}
F^{\ast }(x,\alpha )=\Vert \alpha d\Phi _{\Phi (x)}^{-1}\Vert _{\ast
},\ x\in M,\alpha \in T_{x}^{\ast }M.  \label{polar-minko}
\end{equation}%
If we consider the change of variables $\tilde x=\Phi (x)$, relations (\ref{finsler-derivaltak}) and (\ref%
{polar-minko}) imply
\begin{equation}
|\nabla v|_{d_F}(x)= F^{\ast }(x,Dv(x))=\Vert (D(v\circ \Phi
^{-1})(\tilde x))\Vert_{\ast },\ v\in C_0^\infty(M),\ x\in M.
\label{deriv-isom}
\end{equation}%
Thus,  for every $v\in C_0^\infty(M)$, $p\in (1,n)$ and $q>0$, we
have
\begin{eqnarray*}
 \|D(v\circ \Phi ^{-1})\|_{L^p(\mathbb R^n,d\tilde
x)}&=& \left(\int_{\mathbb R^n}\Vert (D(v\circ \Phi ^{-1})(\tilde
x))\Vert_{\ast }^pd\tilde
x\right)^\frac{1}{p}=\left(\int_{M}(|\nabla
v|_{d_F}(x))^pdV_F(x)\right)^\frac{1}{p} \\
   &=& \||\nabla
v|_{{{d_F}}}\|_{L^p(M,dV_F)},
\end{eqnarray*}
$${\bf Ent}_{d\tilde x}(|v\circ \Phi
^{-1}|^p)={\bf Ent}_{dV_F}(|v|^p) \ {\rm and}\ \|v\circ \Phi
^{-1}\|_{L^q}=\|v\|_{L^q}.$$ It remains to apply the results of
Cordero-Erausquin, Nazaret and Villani \cite{CE-N-Villani} and
Gentil \cite{Gentil} (cf. Theorems A \& B) for $u=v\circ \Phi
^{-1}$. \hfill $\square$\\

\vspace{0.1cm}

\noindent {\bf Acknowledgments.} The author is grateful to the
Universit\"at Bern for the warm hospitality where this work has been
initiated. He thanks Professor Zolt\'an M. Balogh and Professor
C\'edric Villani for stimulating conversations on the
$L^p-$logarithmic Sobolev inequality 
and Professor Michel Ledoux for suggesting the study
of the whole family of Gagliardo-Nirenberg inequalities. Research
supported by a grant of the Romanian National Authority for
Scientific Research, CNCS-UEFISCDI,  "Symmetries in elliptic
problems: Euclidean and non-Euclidean techniques", project no.
PN-II-ID-PCE-2011-3-0241, and by the J\'anos Bolyai Research
Scholarship of the Hungarian Academy of Sciences.


\end{document}